\documentclass[journal]{IEEEtran}
\usepackage
{
    graphicx,
    amssymb,
    amsmath,
    amsthm,
    xcolor,
    dsfont,
    algpseudocode,
    authblk,
    subcaption
}

\usepackage[ngerman, english]{babel}
\usepackage{bm}
\usepackage{enumitem}
\usepackage{algorithm}
\usepackage{tabularx}
\usepackage[colorlinks,
            pdffitwindow=false,
            plainpages=false,
            pdfpagelabels=true,
            pdfpagemode=UseOutlines,
            pdfpagelayout=SinglePage,
            bookmarks=false,
            colorlinks=true,
            hyperfootnotes=false,
            linkcolor=blue,
            citecolor=blue]
{hyperref}

\usepackage{todonotes}
\usepackage{siunitx}

\newcommand{\R}{\mathbb{R}}

\newcommand{\Dcal}{\mathcal{D}}

\newcommand{\Xcal}{\mathcal{X}}
\newcommand{\Ycal}{\mathcal{Y}}

\newcommand{\Pset}{\mathcal{P}}
\newcommand{\Pfront}{\Pset_{\mathcal{F}}}
\newcommand{\Psethat}{\widehat{\Pset}}
\newcommand{\Pfronthat}{\Psethat_{\mathcal{F}}}

\newcommand{\iterate}[2]{{#1}^{(#2)}}

\newcommand{\norm}[1]{\left\|#1\right\|}

\newcommand{\condSet}[2]{\left\{ #1 ~ \middle\vert~ #2 \right\}}

\begin{document}
%
\title{Surrogate-assisted multi-objective design of complex multibody systems}
%
%
\author{Augustina C. Amakor, Manuel B. Berkemeier, Meike Wohlleben, Walter Sextro, Sebastian Peitz
\thanks{ACA, MBB and SP are with the Department of Computer Science, TU Dortmund,
Dortmund, Germany, and with the Lamarr Institute for Machine Learning and Artificial Intelligence, e-mail: \{augustina.amakor,manuel.berkemeier,sebastian.peitz\}@tu-dortmund.de. MW and WS are with the Faculty of Mechanical Engineering, Paderborn University, Paderborn, Germany, e-mail:\{meike.wohlleben,walter.sextro\}@upb.de}
}

%

\maketitle

\begin{abstract}
The optimization of large-scale multibody systems is a numerically challenging task, in particular when considering multiple conflicting criteria at the same time. In this situation, we need to approximate the Pareto set of optimal compromises, which is significantly more expensive than finding a single optimum in single-objective optimization. To prevent large costs, the usage of surrogate models---constructed from a small but informative number of expensive model evaluations---is a very popular and widely studied approach. The central challenge then is to ensure a high quality---that is near-optimality---of the solutions that were obtained using the surrogate model, which can be hard to guarantee with a single pre-computed surrogate. We present a back-and-forth approach between surrogate modeling and multi-objective optimization to improve the quality of the obtained solutions. Using the example of an expensive-to-evaluate multibody system, we compare different strategies regarding multi-objective optimization, sampling and also surrogate modeling, to identify the most promising approach in terms of computational efficiency and solution quality.
\end{abstract}


\section{Introduction}\label{sec:Introduction}
A large number of complex systems is composed of interconnected rigid or flexible bodies that interact with each other through joints, forces, or constraints, also referred to as \emph{multibody systems} \cite{Schiehlen1997,Schiehlen2007}. Examples range from robots over vehicle suspension systems or wind turbines to human biomechanics.
Systems of this type are often modeled mathematically by a set of coupled ordinary (or partial) differential equations constrained by algebraic constraints.
With increasing numbers of components, the numerical simulation of these systems becomes increasingly expensive, rendering multi-query tasks such as optimization or control very costly. This issue is further amplified in the context of \emph{multi-objective optimization} \cite{Miettinen1998,Ehr05}. There, we want to optimize multiple conflicting criteria at the same time, such that instead of a single optimum, we need to calculate the \emph{Pareto set} of optimal compromises, that is, an entire set of optimal points.

As solving \emph{multi-objective optimization problems} (\emph{MOPs}) is computationally expensive---in particular for complex and costly-to-evaluate models---there is a strong interest in accelerating the evaluation of the objective function or their gradients.
There has been extensive research on surrogate-assisted multi-objective optimization \cite{Chugh2015,Tabatabaei2015,PD18a,Deb2021}. That is, instead of the expensive model, we train a surrogate function using a small number of expensive model evaluations at carefully chosen points.
Modeling techniques for range from polynomials over radial basis functions and Kriging models to neural networks, and there is a distinction between global approximations and ones that are valid only locally. Besides smaller models, the latter case may allow for error analysis through trust-region techniques \cite{Berkemeier2021}, at the cost of requiring additional intermittent evaluations of the original objective function.
Not surprisingly, modern machine learning has found its entrance into this area of research as well, see \cite{Qu2021,PH24} for recent overviews. 
Deep surrogate models were suggested in \cite{Botache2024,Zhang2022}, and generative Kriging modeling was studied in \cite{Hussein2016}. Another generative modeling approach called GFlowNets was proposed in \cite{Jain2023}, and the usage of large language models was suggested in \cite{Liu2023}.

In this paper we present an algorithm (in Section \ref{sec:Methodology}) for calculating the entire Pareto front as accurately as possible, while respecting a prescribed budget of expensive calculations, with a specific focus on expensive-to-simulate multibody systems (the basics are introduced Section \ref{subsec:MBS}). We compare 
\begin{itemize}
    \item a single offline training phase with a back-and-forth procedure between surrogate modeling (Section \ref{subsec:Surrogate}) and multi-objective optimization (\ref{subsec:MO}).
\end{itemize}
In the latter case, surrogate-based Pareto sets (using either radial basis functions or deep neural networks as surrogate models) can be used for more informed sampling. Furthermore, we compare  
\begin{itemize}
    \item multi-objective evolutionary algorithms and gradient descent with multi-start, as well as
    \item surrogate modeling via radial basis functions or deep neural networks.
\end{itemize} 
Our findings (Section \ref{sec:Results}) suggest that an interactive coupling between the well-known NSGA-II algorithm \cite{Deb2002} and a moderate number of Pareto-informed sampling steps for neural network surrogates is the best way to find near-optimal solutions.

\section{Preliminaries}\label{sec:Preliminaries}

\subsection{Multibody systems}\label{subsec:MBS}
To analyze the dynamic behavior of a mechanical system, a physical model in the form of mathematical equations is required. These models aim to represent the system under investigation as accurately as necessary but as simply as possible. \\
Multibody systems (MBS) consist of three fundamental elements: bodies with mass, joints connecting these bodies, and forces or moments acting on the bodies. The bodies can be either rigid or flexible, but for simplicity, we will focus on rigid bodies in this context. Likewise, we restrict ourselves here to translations and the conservation of impulse, although the same approach is equally valid for rotations using the conservation of momentum.
Forces in the system are classified as either applied forces, such as those generated by massless springs or dampers, or reaction forces, which represent the internal interactions at the joints. Using the Newton-Euler formalism, the equations of motion for such a system can be formulated as the ordinary differential equation
\begin{equation}\label{eq:MBS}
    M(y,t) \ddot{y}(t)+C(y,\dot{y},t) = Q(y,\dot{y},t).
\end{equation}
Here, $M(y,t)$ represents the mass matrix, $C(y, \dot{y}, t)$ is the vector of generalized Coriolis, centrifugal, and gyroscopic forces (i.e., reaction forces), and $Q(y, \dot{y}, t)$ is the vector of generalized applied forces. The vector $y$ denotes the positions, $\dot{y} = \frac{\mathrm{d}y}{\mathrm{d}t}$ the velocities, and $\ddot{y} = \frac{\mathrm{d}^2y}{\mathrm{d}t^2}$ the accelerations. By solving this differential equation, typically numerically, the dynamic behavior of the described MBS can be simulated and analyzed \cite{Woernle2024,Schiehlen1997}.

\subsection{Multi-objective optimization}\label{subsec:MO}
This section covers only the basic concepts of multi-objective optimization, for detailed introductions, we refer the reader to \cite{Miettinen1998,Ehr05}.
Consider the situation where instead of a single loss function, we have a vector with $K$ conflicting ones, i.e., $f(x) = [f_1(x), \ldots, f_K(x)]^\top$. The task thus becomes to minimize all losses at the same time, i.e.,
\begin{equation}\label{eq:MOP}
    \min_{x\in\R^N} \begin{pmatrix}f_1(x) \\ \vdots \\ f_K(x)\end{pmatrix}. \tag{MOP}
\end{equation}
If the objectives are conflicting, then there does not exist a single optimal $x^*$ that minimizes all $f_k$. Instead, there exists a \emph{Pareto set} $\Pset$ with optimal trade-offs, i.e., 
\[
\resizebox{\linewidth}{!}{$
\Pset = \condSet{x\in\R^N}{\nexists\hat{x}:\begin{array}{ll}
    f_k(\hat{x}) \leq f_k(x) & \mbox{for}~k=1,\ldots,K, \\
    f_k(\hat{x}) < f_k(x) & \mbox{for at least one }~k
\end{array}}.
$}
\]
In other words, a point $\hat{x}$ \emph{dominates} a point $x$, if it is at least as good in all $f_k$, while being strictly better with respect to at least one loss. The Pareto set $\Pset$ thus consists of all non-dominated points.
The corresponding set in objective space is called the \emph{Pareto front} $\Pfront=f(\Pset)$. Under smoothness assumptions, both objects have dimension $K-1$ \cite{Hil01}, i.e., $\Pset$ and $\Pfront$ are lines for two objectives, 2D surfaces for three objectives, and so on. Furthermore, they are bounded by Pareto sets and fronts of the next lower number of objectives \cite{GPD19}, meaning that individual minima constrain a two-objective solution, 1D fronts constrain the 2D surface of a $K=3$ problem, etc.

Closely related to single-objective optimization, there exist first order optimality conditions, referred to as the \emph{Karush-Kuhn-Tucker (KKT)} conditions \cite{Miettinen1998}. A point $x^*$ is said to be \emph{Pareto-critical} if there exists a convex combination of the individual gradients $\nabla f_k$ that is zero. More formally, we have
\begin{align}
    \sum_{k=1}^K \alpha^*_k \nabla f_k(x^*) =0, \qquad \sum_{k=1}^K \alpha^*_k =1, \label{eq:KKT} \tag{KKT}
\end{align}
which is a natural extension of the case $K=1$. 

Equation \eqref{eq:KKT} is the basis for most gradient-based methods. Their goal is to compute elements from the Pareto critical set,
\[
    \Pset_c = \condSet{x^*\in\R^N}{\exists \alpha^*\in\R_{\geq 0}^{K},\, \sum_{k=1}^K \alpha^*_k =1\,:~ \eqref{eq:KKT}\,\mbox{holds}},
\]
which contains excellent candidates for Pareto optima, since $\Pset_c \supseteq \Pset$. In a similar fashion to the single-objective case, there exist extensions to constraints \cite{GPD17,Fliege2016}, but we will exclusively consider unconstrained problems here.

\subsubsection*{Overview of methods}
There exists a large number of conceptually very different methods to find (approximate) solutions of \eqref{eq:MOP}. The first key question is whether we put the decision-making before or after optimization, or even consider an interactive approach.
To calculate single points, \emph{multiple-gradient descent algorithms (MGDAs)} are increasingly popular, in particular when it comes to very high-dimensional problems. The key ingredient is the calculation of a \emph{common descent direction} $d(x)\in\R^N$ that satisfies
\[
    \left(\nabla f_k(x)\right)^\top d(x) < 0, \quad k\in\{1,\ldots,K\},
\]
which again is a straightforward extension of single-objective descent directions. The determination of such a $d$ usually requires the solution of a subproblem in each step, for instance a quadratic problem of dimension $K$ \cite{Schaeffler2002,Dsidri2012},
\begin{equation}\tag{CDD}\label{eq:CDD}
\begin{aligned}
    d(x) &= -\sum_{k=1}^K w_k \nabla f_k(x),\qquad \mbox{where}  \\
    w&=\arg\min_{\begin{array}{c}\scriptstyle{\hat{w}\in[0,1]^K} \\ \scriptstyle{\sum \hat w_k=1} \end{array}} \left\|\sum_{k=1}^K \hat w_k \nabla f_k(x)\right\|_2^2.
\end{aligned}
\end{equation}
However, there are various alternatives to \eqref{eq:CDD} such as a dual formulation \cite{FS00} or a computationally cheaper so-called \emph{Franke-Wolfe} approach \cite{SK18}. Once a common descent direction $d(x)$ has been obtained, we proceed in a standard fashion by iteratively updating $x$ until convergence or some other stopping criterion is met, cf.\ Algorithm \ref{alg:MGDA}.
\begin{algorithm}[thb]
\caption{Multiple-gradient descent algorithm (MGDA)}\label{alg:MGDA}
\begin{algorithmic}[1]
\Require Initial guess $x^{(0)}$, learning rate $\eta \in\R_{>0}$ (possibly adaptive), maximum number of iterations $i_{\mathsf{max}}$, hyperparameters (depending on specific version of MGDA)
\Ensure $x^*\in\Pset_c$
\State Set $i=0$
\While{$\iterate{x}{i}\notin\Pset_c$ \textbf{and} $i<i_{\mathsf{max}}$}
    \State Calculate gradients $\nabla f_i\left(\iterate{x}{i}\right)$ for $i=1\ldots,k$
    \State Calculate descent direction $d\left(\iterate{x}{i}\right)$ (e.g., via \eqref{eq:CDD})
    \State If adaptive, determine learning rate $\eta\left(\iterate{x}{i}\right)$
    \State Update $x$:
    \[
        \iterate{x}{i+1}= \iterate{x}{i} + \eta\left(\iterate{x}{i}\right) d\left(\iterate{x}{i}\right)
    \]
    \State $i = i + 1$
\EndWhile
\end{algorithmic}
\end{algorithm}
Various extensions concern Newton \cite{Fliege2009} or Quasi-Newton \cite{Povalej2014} directions, uncertainties \cite{PD18b}, momentum \cite{Tanabe2023,SP24a,Nikbakhtsarvestani2023}, or non-smoothness \cite{Miettinen1995,GP21,Tanabe2018,Tanabe2023}. 

\begin{figure}[thb]
    \centering
    \includegraphics[width=.9\linewidth]{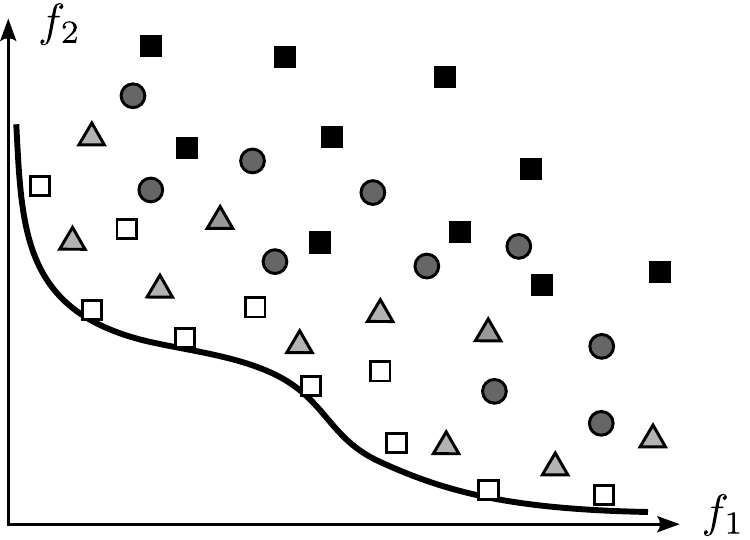}
    \caption{MOEA example, where a population of individuals is improved from one generation to the next ($\blacksquare\rightarrow\bigcirc\rightarrow\triangle\rightarrow\square$).}
    \label{fig:MOEA}
\end{figure}

If we want to postpone the decision-making to take a more informed decision, we need to calculate the entire Pareto set $\Pset$ and front $\Pfront$. The most straightforward approach is to adapt the weights $w$ in scalarization and solve the single-objective problem multiple times. Alternatively, one can combine MGDA with a multi-start strategy (i.e., a set of random initial guesses $\left\{x^{(0,j)}\right\}_{j=1}^M$) to obtain multiple points. However, in both cases, it may be very hard or even impossible to obtain a good coverage of $\Pset$, i.e., that approximates the entire set with evenly distributed points.
Instead, we can directly consider a \emph{population} of weights $\left\{x^{(j)}\right\}_{j=1}^M$ that we iteratively update to improve each individual's performance while also ensuring a suitable spread over the entire front $\Pfront$. 
Algorithms of this class are referred to as \emph{multi-objective evolutionary algorithms (MOEAs)} \cite{CLV07}, see Figure \ref{fig:MOEA} for an illustration and Algorithm \ref{alg:MOEA} for a rough algorithmic outline. Therein, the population's fitness is increased from generation to generation by maximizing a criterion that combines optimality (or non-dominance) with a spreading criterion. The most popular and widely used algorithm in this category is likely NSGA-II \cite{Deb2002}, but there are many alternatives regarding the crossover step (Step 3 in Algorithm \ref{alg:MOEA}), mutation (Step 4), the selection (Step 5), the population size, and so on. For more details, see the surveys \cite{Fonseca1995,Zhou2011,Tian2021}. Finally, many combinations exist with---for instance---preference vectors \cite{Thiele2009} or gradients \cite{Bosman2012,Nikbakhtsarvestani2023}, such that the creation of offspring is more directed using gradient information.

\begin{algorithm}[thb]
\caption{Multi-objective evolutionary algorithm (MOEA)}\label{alg:MOEA}
\begin{algorithmic}[1]
\Require Initial population $P(0)$ of individuals $\left\{x^{(0,j)}\right\}_{j=1}^M$, number of generations $i_{\mathsf{max}}$,  hyperparameters (depending on specific version of MOEA)
\State Set $i=0$
\While{$i<i_{\mathsf{max}}$}
    \State Create an offspring population $\widehat{P}{(i)}$ out of ${P}{(i)}$, e.g., using \textbf{crossover} between two individuals
    \State Modify offspring population via \textbf{mutation}:
    \[\widetilde{P}{(i)} = \mathcal{M}\left( \widehat{P}{(i)} \right)\]
    \State \textbf{Selection} of the next generation ${P}{(i+1)}$ either from $\widetilde{P}{(i)}$ or from ${P}{(i)} \cup \widetilde{P}{(i)}$ (the latter is called \emph{elitism}) by a survival-of-the-fittest process (e.g., using a non-dominance and spread metric)
    \State $i = i + 1$
\EndWhile
\end{algorithmic}
\end{algorithm}

\subsection{Surrogate modeling}\label{subsec:Surrogate}
Solving MOPs is naturally more expensive than identifying a minimizer of a single-objective problem. Consequently, expensive models $f$ quickly result in prohibitively large computational cost. It is thus a straightforward idea to replace $f: \R^N \rightarrow \R^K$ by a surrogate model $g: \R^N \times \R^P \rightarrow \R^K$ which depends on $P$ trainable parameters $\theta\in\R^P$. 
There is a lot of research on surrogate modeling, either directly for the objective function $f$ (in the multi-objective context see, e.g., \cite{Chugh2015,Tabatabaei2015,Hussein2016,Deb2021,Berkemeier2021,Zhang2022,Jain2023,Liu2023,Botache2024,PH24}), but also for the simulation of complex dynamical systems \cite{BGW15,PD18a}.\footnote{In the recent reinforcement learning literature, such models are also referred to as \emph{world models} \cite{ha2018world,Wu2023}.} 
Given a small training data set $\Dcal=\left\{\left(\iterate{x}{i},f\left(\iterate{x}{i}\right)\right)\right\}_{i=1}^s$ of size $s$, fitting $g$ is a standard supervised learning problem:
\begin{equation}\label{eq:modelFit}
    \min_{\theta\in\R^P} \frac{1}{s} \sum_{i=1}^s \norm{f\left(\iterate{x}{i}\right) - g\left(\iterate{x}{i}; \theta\right)}_2^2.
\end{equation}
Instead of the mean squared loss, one can use various alternatives and also consider additional regularization terms such as the $\ell_1$ or the $\ell_2$ norm.
As a result, we can replace the expensive \eqref{eq:MOP} by a surrogate version
\begin{equation}\label{eq:MOPs}
    \min_{x\in\R^N} \begin{pmatrix}g_1(x; \theta) \\ \vdots \\ g_K(x; \theta)\end{pmatrix}. \tag{$\widehat{\mbox{MOP}}$}
\end{equation}
We will denote the Pareto set and front of \eqref{eq:MOPs} by $\Psethat$ and $\Pfronthat$, and the central goal is to obtain $\Psethat\approx\Pset$ (or alternatively a set of solutions such that $f(\Psethat)\approx f(\Pset)$).
A general theme is that the closer we want to get to the true set, the more costly the surrogate modeling, meaning that this is in itself a multi-criteria problem.
Aside from the specific formulation of the loss function in \eqref{eq:modelFit}, central questions that need to be addressed in this context are:
\begin{itemize}
    \item Which type of surrogate model should we use?
    \item What is the appropriate choice for $s$? (Small $s$ results in low computational effort, at the cost of poor generalization towards points $(x,f(x)) \notin \Dcal$.)
    \item Where should we sample?
    \item How can we ensure convergence or optimality?
\end{itemize}
In the following section, we will address all these points. We will experiment with various values for $m$. In particular, we propose an adaptive back-and-forth strategy between surrogate modeling and multi-objective optimization that helps us to choose samples in a more informed way, while ensuring convergence. We will test our method using both MGDA and the popular evolutionary algorithm NSGA-II, as well as surrogate models based on radial basis functions and deep neural networks, respectively.

\section{Methodology}\label{sec:Methodology}

Optimizing complex multibody systems is often very expensive or even infeasible due to high computational costs. As outlined in the previous section, we thus replace $f(x)$ by a cheaper-to-evaluate approximation $g(x, \theta)$ which is trained in a supervised manner according to \eqref{eq:modelFit}.
In this work, we will focus on artificial neural networks (ANNs) and radial basis functions (RBFs) to approximate $g$. 

\begin{figure}[bht]
    \centering
    \includegraphics[width=0.9\linewidth]{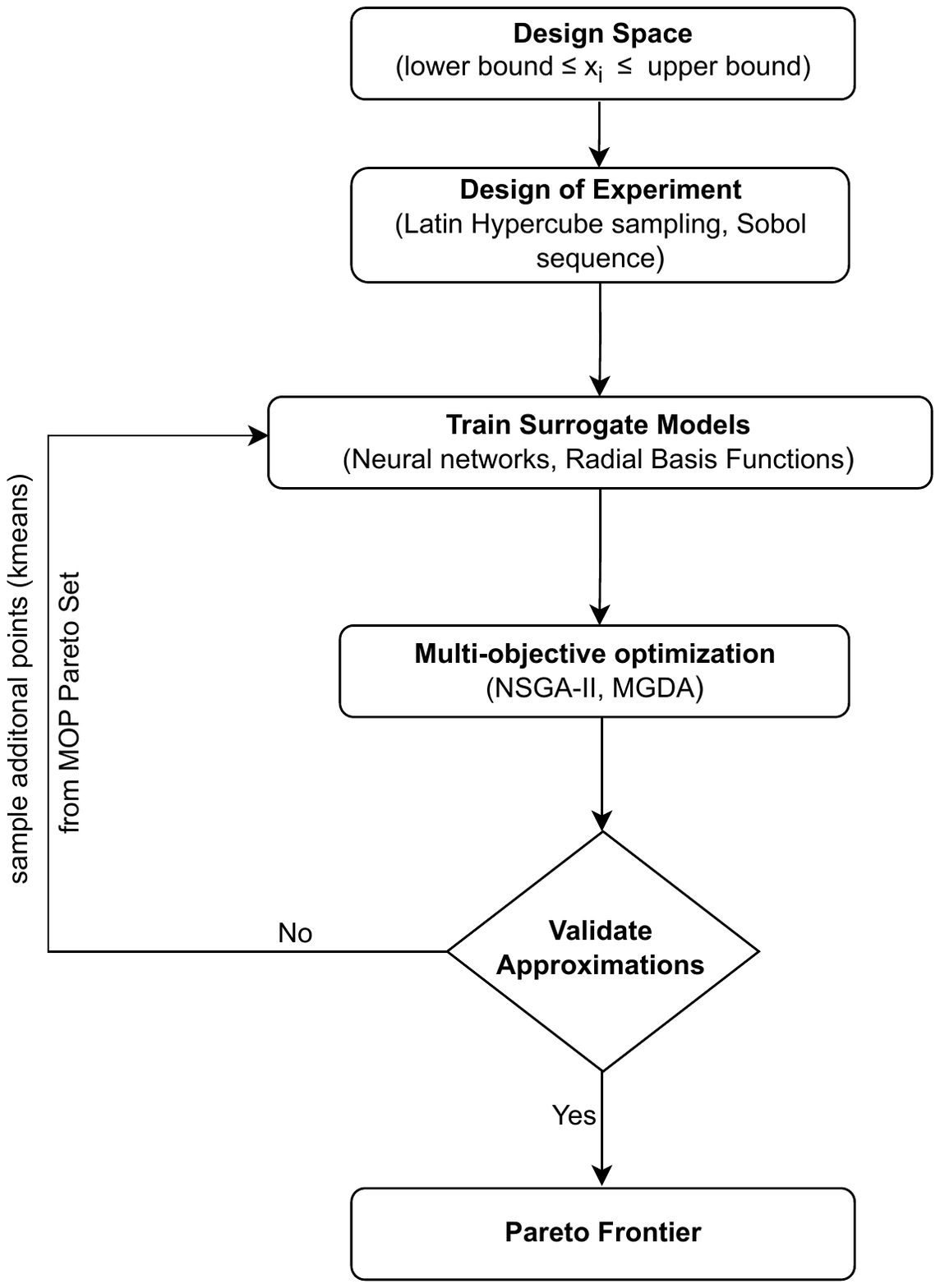}
    \caption{Sketch of the methodology.}
    \label{fig:methodology}
\end{figure}

\begin{figure}[bht]
    \centering
    \includegraphics[width=.8\linewidth]{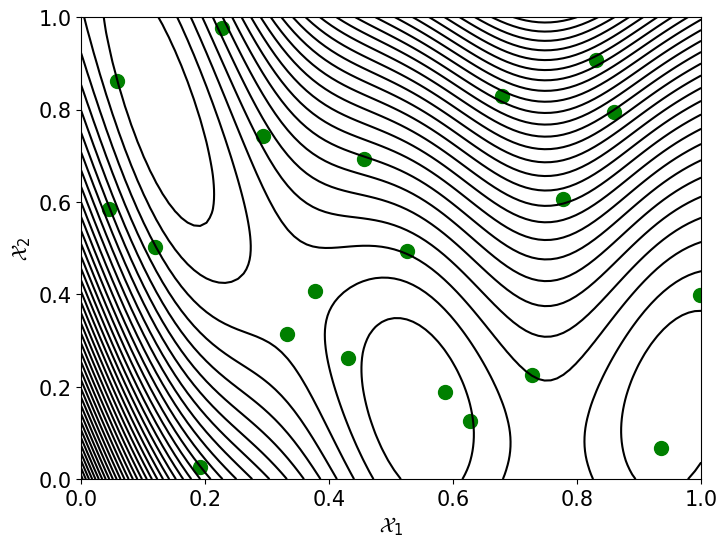}
    \caption{Latin hypercube sampling for 20 points (in green) for the Branin function.}
    \label{fig:lhs}
\end{figure}

The overall procedure is visualized in Figure \ref{fig:methodology} and formalized in Algorithm \ref{alg:SAMO}. The central idea is to intertwine surrogate modeling with optimization, as is very common for trust-region methods \cite{Berkemeier2021}.
In a first step, i.e., where nothing is known about the solution, we collect $s$ samples from the expensive model in a random fashion. Sampling can be performed using various techniques from Design of Experiments (DoE), such as Latin hypercube, Monte-Carlo, Sobol sequences or Halton sets \cite{Thombre2015}. In our case, we are using Latin hypercube sampling, an example of which is presented for the Branin function in Figure \ref{fig:lhs}.



%
\begin{figure}[b!]
    \centering
    \includegraphics[width=.5\linewidth]{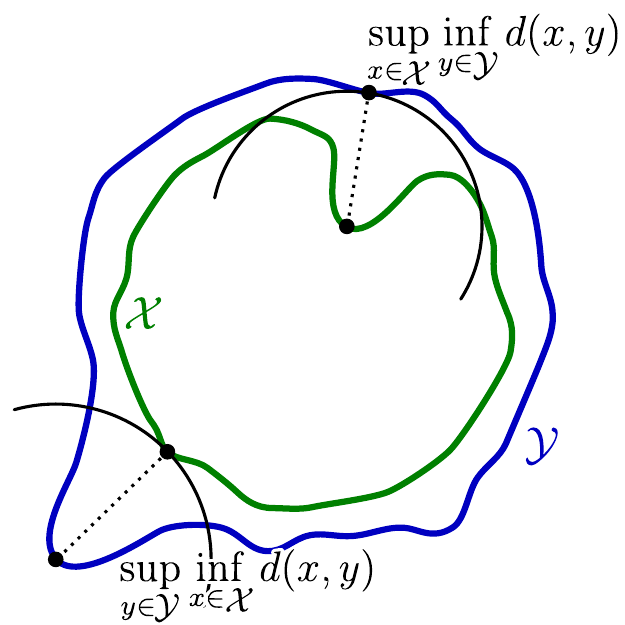}
    \caption{Example for the Hausdorff distance (adapted from \url{https://en.wikipedia.org/wiki/Hausdorff_distance}).}
    \label{fig:Hausdorff}
\end{figure}

Having calculated the first $s$ samples and trained the first surrogate model by solving \eqref{eq:modelFit}, we are already in the position to compute an approximation of $\Pset$ and $\Pfront$ very efficiently by solving \eqref{eq:MOPs} using a suitable algorithm. In our example, this accelerates the problem from approximately 30 seconds for a single function evaluation---which would render optimization infeasible---to just a few seconds for solving \eqref{eq:MOPs} using NSGA-II.
The central issue with this approach is that for a small number $s$, the solution $\Psethat$ and $\Pfronthat$ is likely far away from the true solution $\Pset$ and $\Pfront$. To improve this, we can now compute additional samples, but in a much more informed fashion. Instead of random sampling, we sample $s$ additional samples, distributed equidistantly over the obtained surrogate solution $\Psethat$. For this, we use $k$-means clustering. This now-enhanced data set allows us to find a better approximation of $f$ by $g$, and we can repeat this process until either our budget $S$ of expensive simulations is used up, or until we reach convergence. The latter is realized by measuring the Hausdorff distance (cf.\ \cite{Schuetze2012}) between two sets
\begin{equation}\label{eq:Hausdorff}
\begin{aligned}
    d_H\left( \Xcal, \Ycal \right) = \max \left\{\sup_{y\in \Ycal}\hat{d}(y,\Xcal) , \sup_{x\in \Xcal}\hat{d}(x,\Ycal) \right\}. 
\end{aligned}
\end{equation}
Here, $\hat{d}(y,\Xcal)$ is the Euclidean distance between a point $y\in\Ycal$ and the closest element of the set $\Xcal$, i.e., 
\[\hat{d}(y,\Xcal)=\inf_{x\in\Xcal}d(y,x),\]
where $d$ is the standard Euclidean distance. An equivalent, point-wise definition is visualized in Figure \ref{fig:Hausdorff}. Since we use the surrogate-based Pareto fronts $\iterate{\Pfronthat}{j-1}$ and $\iterate{\Pfronthat}{j}$ for $\Xcal$ and $\Ycal$ in Algorithm \ref{alg:SAMO}---which are finite sets of dimension $M$---the $\sup$ and $\inf$ operators can be replaced by $\max$ and $\min$, respectively.

\begin{algorithm}[t!]
\caption{Surrogate-based adaptive multi-objective optimization algorithm}\label{alg:SAMO}
\begin{algorithmic}[1]
\Require Maximum number of evaluations of the complex model: $S$, Number of samples per iteration: $s\leq S$, Hausdorff distance stopping criterion $h_{\mathsf{min}}$, MOP solver with population size $M$
\Statex \underline{\textit{\textbf{Round 1}: random sampling}}
\State \textit{Sampling} of $s$ random points (e.g., using Latin hypercube) 
\begin{flalign*}\Rightarrow \quad \Dcal=\left\{\left(\iterate{x}{i},f\left(\iterate{x}{i}\right)\right)\right\}_{i=1}^s && \end{flalign*}
\State \textit{Surrogate modeling}: Solve \eqref{eq:modelFit} using $\Dcal$ $\Rightarrow$ $\iterate{g}{0}(x; \theta)$
\State \textit{Optimization}: Solve \eqref{eq:MOPs} using $\iterate{g}{0}(x; \theta)$
\begin{flalign*}\Rightarrow \quad \iterate{\Psethat}{0} \quad\mbox{and}\quad\iterate{\Pfronthat}{0} \quad\mbox{with $M$ elements} && \end{flalign*}
\Statex \underline{\textit{\textbf{Rounds 2, 3, $\ldots$}: Pareto-informed sampling}}
\State Set loop counter $j=1$
\While{$(j-1) \cdot s < S$}
    \State \textit{Sampling}: Determine $s$ sample locations via \emph{$k$-means clustering} of $\iterate{\Psethat}{j-1}$ and collect $s$ new samples in $\widehat{\Dcal}$
    \begin{flalign*}\Rightarrow \quad \Dcal=\Dcal \cup \widehat{\Dcal} && \end{flalign*}
    \State \textit{Surrogate modeling}: Solve \eqref{eq:modelFit} using $\Dcal$ $\Rightarrow$ $\iterate{g}{j}(x; \theta)$
    \State \textit{Optimization}: Solve \eqref{eq:MOPs} $\Rightarrow$ $\iterate{\Psethat}{j}$ and $\iterate{\Pfronthat}{j}$
    \State \textit{Convergence:} Compute \emph{Hausdorff distance} (cf.\ \eqref{eq:Hausdorff})
    \[
    h = d_H\left( \iterate{\Pfronthat}{j-1}, \iterate{\Pfronthat}{j} \right) 
    \]
    \If{$h < h_{\mathsf{min}}$}
    \State \textbf{STOP}
    \EndIf
    \State $j = j + 1$
\EndWhile
\State Calculate the Pareto subsets $\Pset$ and $\Pfront$ using a non-dominance test on $\Dcal$
\end{algorithmic}
\end{algorithm}

As the solution, we can either present the final surrogate-based $\Psethat$ and $\Pfronthat$, or perform a non-dominance test over all real samples (line 15 of Algorithm \ref{alg:SAMO}).


\section{Experiment}\label{sec:Experiment}

\begin{figure}[b!]
    \centering
    \includegraphics[width=.5\linewidth]{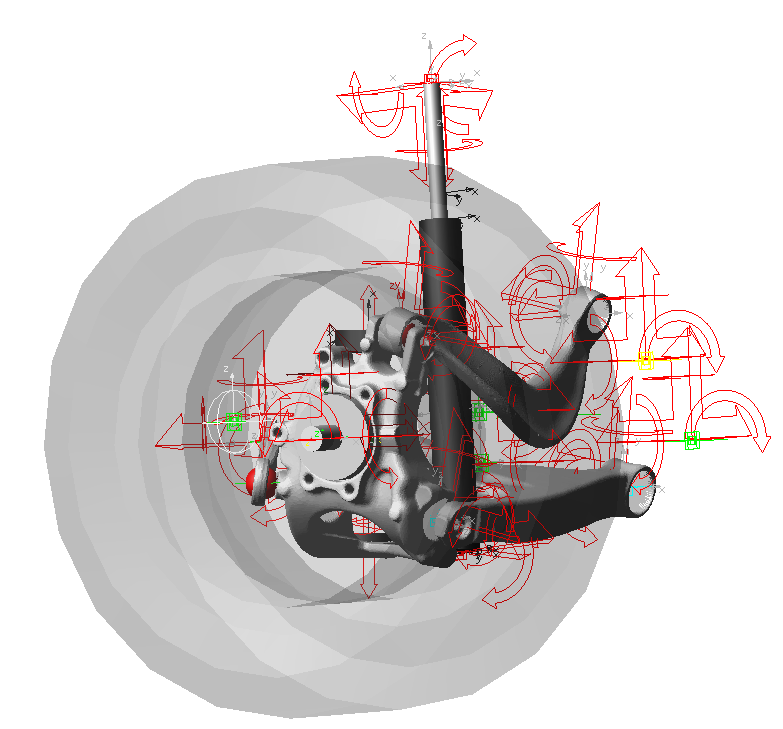}
    \caption{Trapezoidal link rear suspension system under consideration.}
    \label{fig:TLA}
\end{figure}

The MBS we consider as an example is an independent rear suspension system of a passenger car, described in detail in \cite{Schuette2021}, see Figure \ref{fig:TLA}. The wheel suspension consists of a spring, a shock absorber and links connecting the wheel carrier to the subframe of the axle, enabling the transfer of forces and ultimately supporting the vehicle structure. The suspension primarily provides the wheel with a vertical degree of freedom, decoupling the body from, e.g., road excitations, while ensuring a robust contact between tire and road. This is integral to ensuring safe driving behavior, as a car can only be driven safely when the wheel maintains road contact. At the same time, to improve passenger comfort, movements of the vehicle body should remain moderate. The suspension is also crucial for vehicle dynamics, as it defines the wheel's position relative to the road. \cite{heissing2011, Mitschke2014}.

The MBS model was implemented using the Neweul-M$^2$ software, which generates a system of differential equations as seen in Equation \eqref{eq:MBS} \cite{Kurz2010a}. In addition to the elements already mentioned, the model includes a test rig that imposes motion on the wheel hub. The overall model consists of 13 rigid bodies (subframe, chassis, toe link, camber link, trapezoidal link, knuckle, wheel bearing, wheel hub, damper rod, damper housing, lower spring seat, upper spring seat, test rig), 3 point masses (body mass, unsprung mass, brake caliper), 6 force elements (2 spring seats, body spring, damper, rebound stop, test rig), and 22 ideal joints, resulting in a system with 5 remaining degrees of freedom.

For optimization, the joint centers of the links can be translated by a radius of \SI{0.003}{\meter} around their initial position. Toe link and camber link each have two joints, trapezoidal link has four joints, so there are a total of $N = 24$ optimization parameters: $8$ (joints) $\times$ $3$ (spatial directions).
As described, the suspension has a significant impact on various criteria related to driving safety, vehicle dynamics, and ride comfort. To investigate these criteria, the test rig can impose different displacement excitations. However, since only a single wheel suspension is considered, the choice of possible driving maneuvers is limited: lateral dynamics cannot be analyzed as no lateral forces can be absorbed, which precludes cornering maneuvers. Similarly, longitudinal dynamics such as braking and driving cannot be evaluated. We thus focus on the vertical dynamics. Two conflicting criteria are considered: a safety criterion in terms of minimizing wheel load fluctuations, and a comfort criterion in the form of low body accelerations \cite{Mitschke2014,heissing2011}. The driving maneuver is a sinusoidal excitation with an amplitude of $\SI{0.001}{\meter}$ and a frequency of $\SI{7}{\hertz}$ \cite{Mitschke2014}.
The wheel load $F_z$ refers to the vertical forces acting on the tire (in our case, on the wheel hub), and is therefore an entry of the vector $C$ (cf.\ Equation \eqref{eq:MBS}), which contains the reaction forces. Wheel load fluctuations are the variations in wheel load at a single wheel caused by driving on uneven roads \cite{Mitschke2014}. 
Therefore, the safety objective function is the amplitude of $F_z$, estimated over the simulated time horizon $t_0$ to $t_e$ by 
\begin{equation}
    \hat{F}_z = \frac{1}{2}\left(\max_{t\in[t_0,t_e]}F_z(t)-\min_{t\in[t_0,t_e]}F_z(t)\right).
    \label{equ:F_z}
\end{equation} 
The body acceleration $\ddot{z}_A$ represents the vertical acceleration of the vehicle body, and is an entry of the acceleration vector $\ddot{y}$. The amplitude
\begin{equation}
    \hat{\ddot{z}}_A = \frac{1}{2} \left(\max_{t\in[t_0,t_e]}\ddot{z}_A(t)-\min_{t\in[t_0,t_e]}\ddot{z}_A(t)\right)
    \label{equ:Z_A}
\end{equation}
is selected as the comfort objective function \cite{Mitschke2014}. We thus have the following multi-objective optimization problem of the form \eqref{eq:MOP}.
\begin{equation}\label{eq:MOP_MBS}
\min_{x\in\R^{24}}f(x) = \min_{x\in\R^{24}}\begin{pmatrix}\hat{F}_z\\\hat{\ddot{z}}_A\end{pmatrix}.
\end{equation}
The objective function $f$ is evaluated in a black-box fashion, i.e., we set the values for $x$, then perform a simulation over the time horizon $[t_0,t_e]$, collect (discretized) trajectories $F_z(t)$ and $\ddot{z}_A(t)$, and finally calculate $f_1(x)$ and $f_2(x)$.

Given the expensive nature of optimizing the two objectives of the MBS expressed in equations \eqref{equ:F_z} and \eqref{equ:Z_A}, the surrogate-based adaptive multi-objective optimization algorithm in Algorithm \ref{alg:SAMO} is used. The experiments are carried out on a machine with 2.10 GHz 12th Gen Intel(R) Core(TM) i7-1260P CPU and 32 GB memory, using Python 3.11.1 and MATLAB R2024b.
Our implementation is in Python, but calls MATLAB to evaluate \eqref{eq:MOP_MBS} in a black-box fashion.
Further settings are:
\begin{itemize}
    \item the design space for the optimization parameters of MBS is set to the infinity norm i.e., $-0.003\le x_i \leq0.003$.
    \item as described in Algorithm \ref{alg:SAMO} the Latin hypercube is used for initial sampling and the $k$-means clustering for succeeding sampling.  
    \item we found $h_{\mathsf{min}}=2$ to be a good value for our particular application.
\end{itemize}

\section{Results}\label{sec:Results}
We now investigate the performance of Algorithm \ref{alg:SAMO} for the efficient surrogate-based solution of Problem \eqref{eq:MOP_MBS}. First, we show a detailed analysis over multiple iterations using NSGA-II and ANN surrogate modeling. We then study different variations in terms of surrogate modeling and sample sizes $s$. In our experiments, we have also investigated MGDA (Algorithm \ref{alg:MGDA}) with multi-start as an alternative to NSGA-II, but we found it to be inferior both in terms of the spread across the Pareto front as well as in terms of computing times. We thus only consider NSGA-II as the optimizer from now on.

\subsection{NSGA-II optimizer with ANN surrogate models}

\begin{figure}
    \centering
    \includegraphics[width=0.8\linewidth]{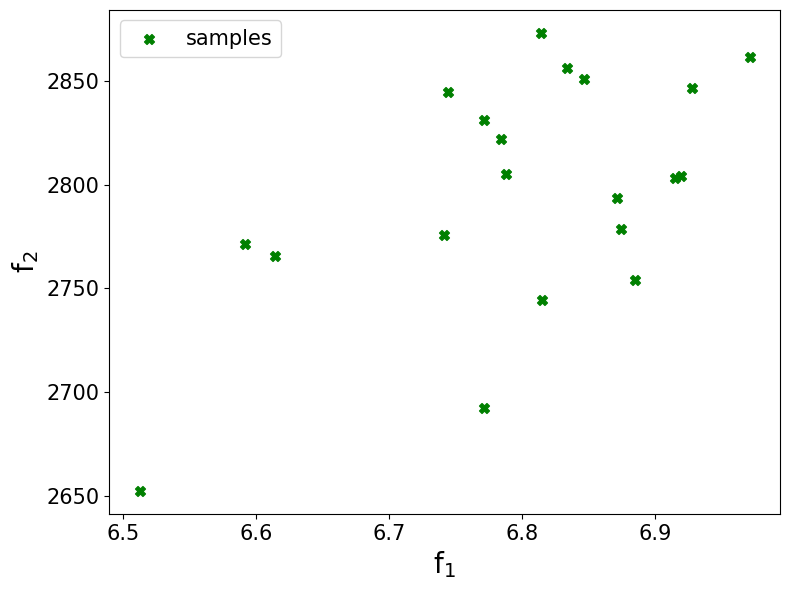}
    \caption{The first $s=20$ points computed by evaluating the complex model \eqref{eq:MOP_MBS} using the Latin hypercube sampling.}
    \label{fig:lhs20}
\end{figure}
\begin{figure}
    \centering
    \includegraphics[width=0.8\linewidth]{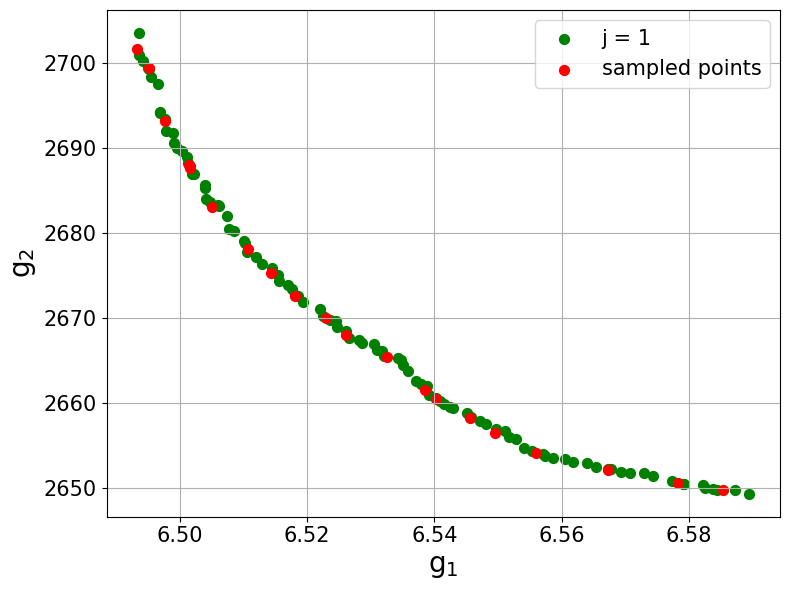}
    \caption{The Pareto front $\Pfronthat$ using the ANN surrogate model (constructed using the samples in Figure \ref{fig:lhs20}) and NSGA-II are shown in green color. In red are the first $s=20$ Pareto-informed samples, obtained via $k$-means clustering of $\Pfronthat$.}
    \label{fig:sampledKM}
\end{figure}

To analyze the performance and demonstrate the efficiency of our approach, we first use a fully connected neural network with two hidden layers of 64 neurons each 
as our surrogate model. For now, we set the sample size to $s=20$. In the first loop of Algorithm \ref{alg:SAMO}, we begin by sampling $s$ times from $-0.003\le x_i \leq0.003$ using Latin hypercube, and compute the true values $f(x)$, cf.\ Figure \ref{fig:lhs20}. These 20 points are scaled (i.e., normalized) and split in a $80:20$ ratio for validation. We then train the ANN surrogate model using the ADAM optimizer \cite{Adam2015} and the mean squared error (cf.\ Equation \eqref{eq:modelFit}) as the error metric. 
We solve the surrogate-based problem \eqref{eq:MOPs} using NSGA-II with a population size of $M=100$, $i_{\mathsf{max}}=200$ generations, a simulated binary crossover of probability $0.5$ and polynomial mutation with distribution index $\eta=20$. For the sampling of new data points, we use $k$-means clustering. Figure \ref{fig:sampledKM} shows the first approximate Pareto front $\Pfronthat$ and the location of the $20$ Pareto-informed samples.

\begin{figure*}[h!]
    \begin{subfigure}[b]{0.5\columnwidth}
        \includegraphics[width=\linewidth]{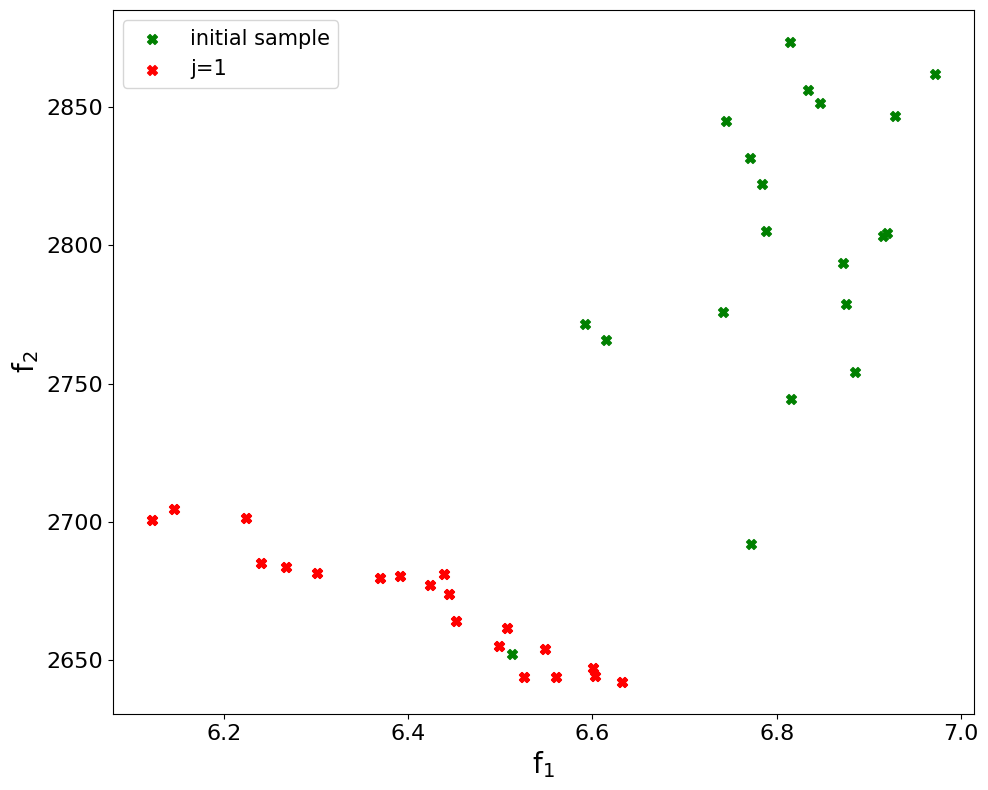}
        \caption{}
        \label{fig:sub1_n20}
    \end{subfigure}%
    \begin{subfigure}[b]{0.5\columnwidth}
        \includegraphics[width=\linewidth]{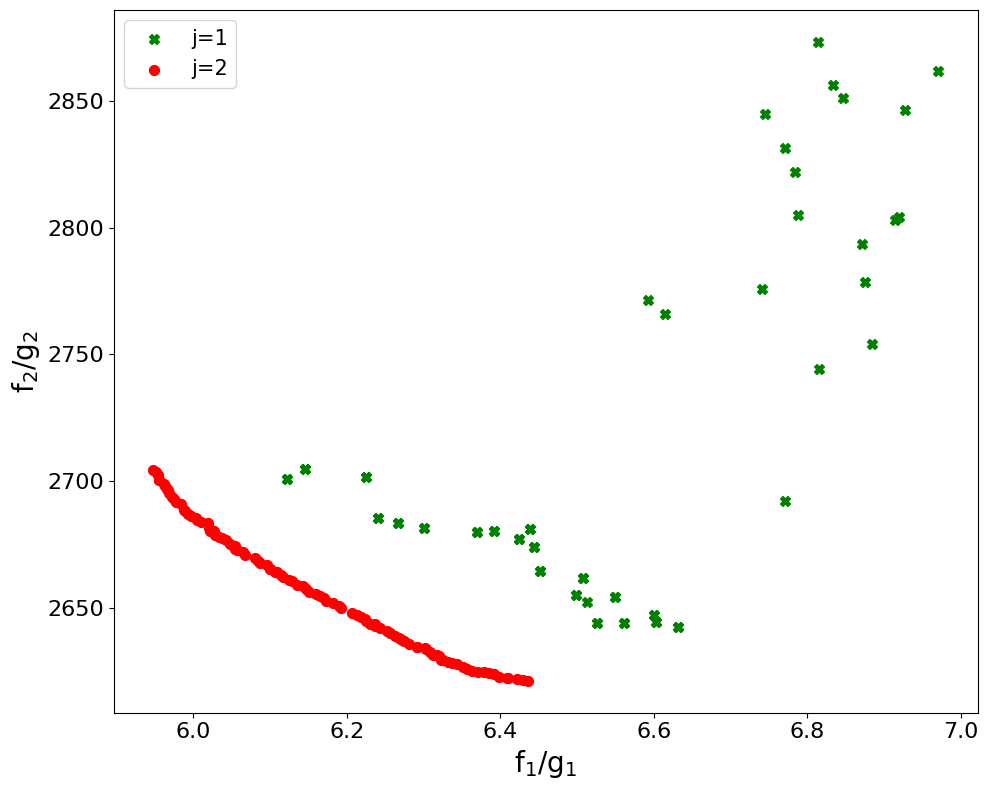}
        \caption{}
        \label{fig:sub2_n20}
    \end{subfigure}%
    \begin{subfigure}[b]{0.5\columnwidth}
        \includegraphics[width=\linewidth]{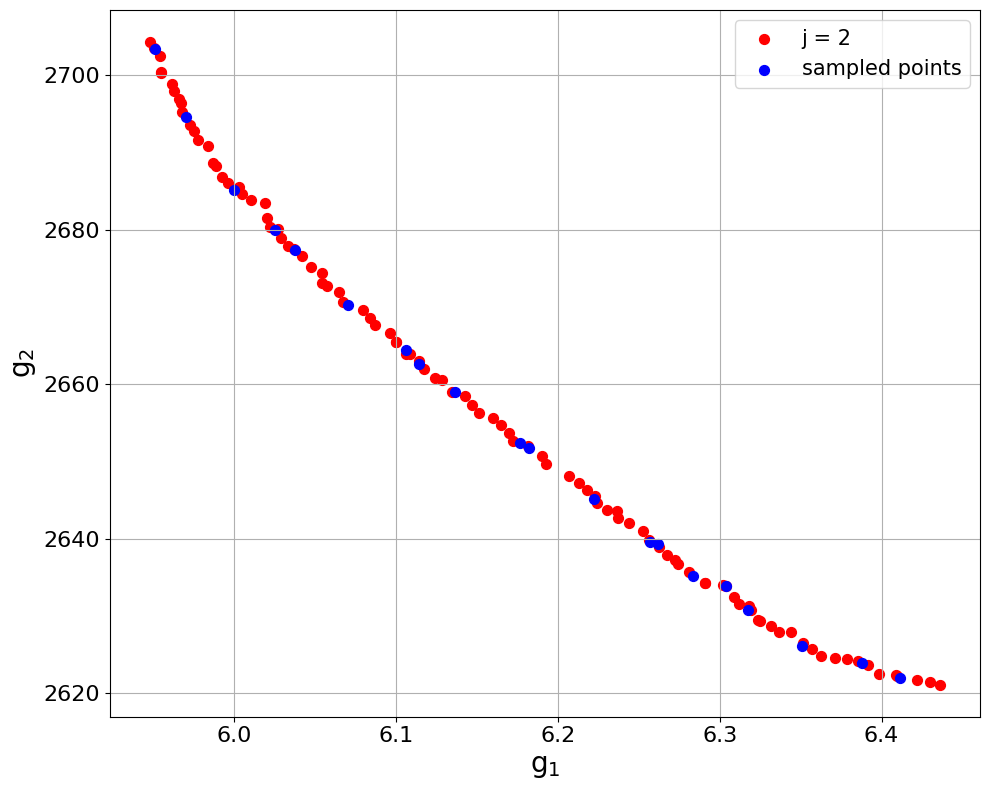}
        \caption{}
        \label{fig:sub3_n20}
    \end{subfigure}%
    \begin{subfigure}[b]{0.5\columnwidth}
        \includegraphics[width=\linewidth]{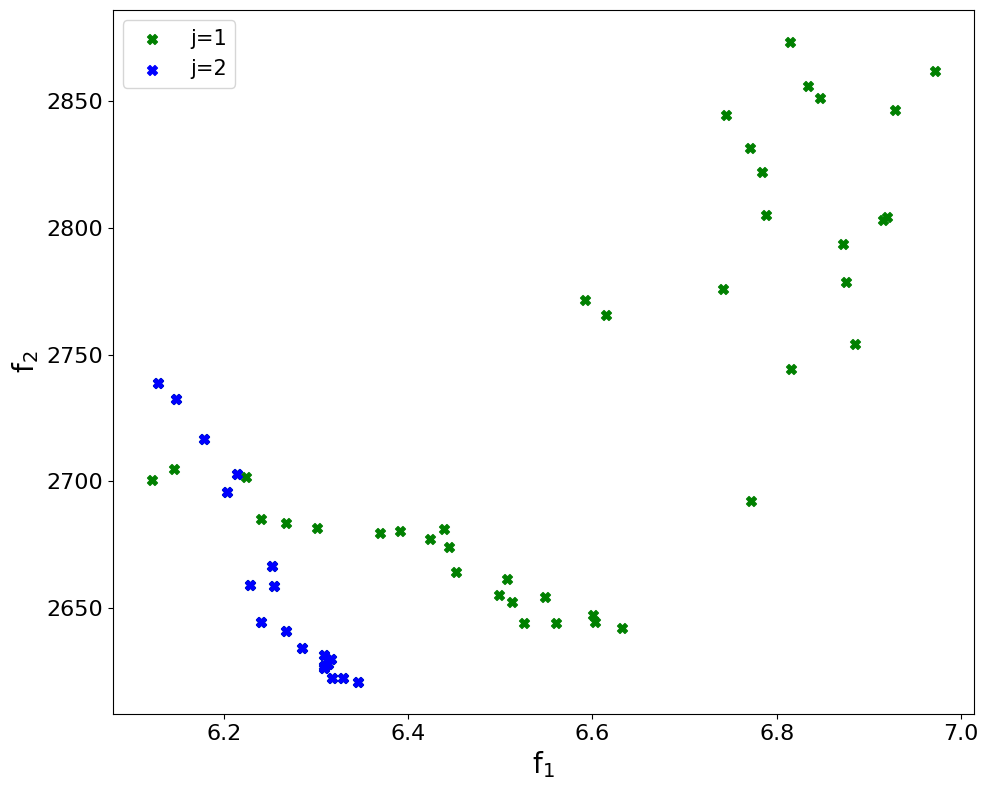}
        \caption{}
        \label{fig:sub4_n20}
    \end{subfigure}
    \begin{subfigure}[b]{0.5\columnwidth}
        \includegraphics[width=\linewidth]{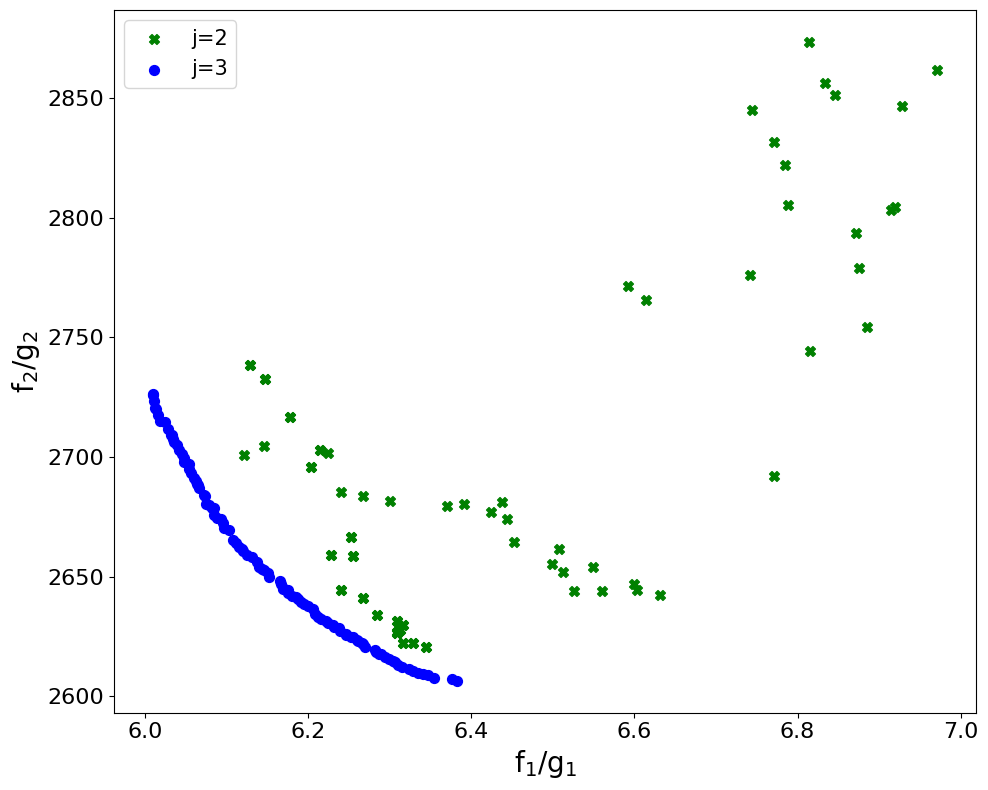}
        \caption{}
        \label{fig:sub5_n20}
    \end{subfigure}%
    \begin{subfigure}[b]{0.5\columnwidth}
        \includegraphics[width=\linewidth]{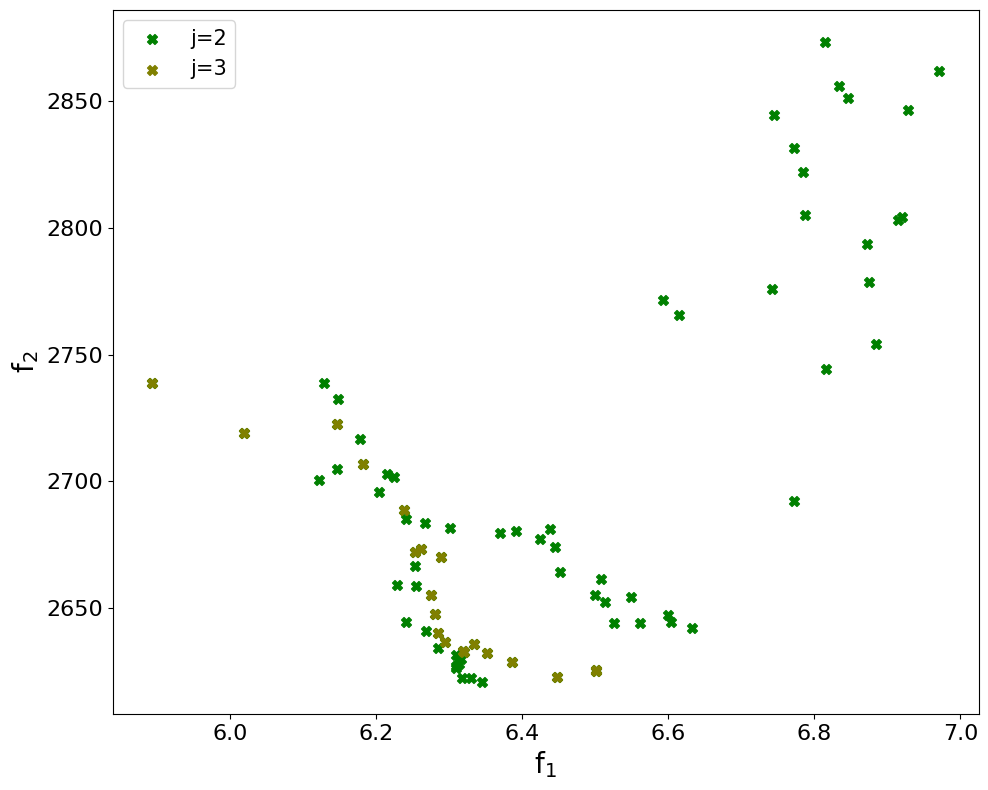}
        \caption{}
        \label{fig:sub6_n20}
    \end{subfigure}%
    \begin{subfigure}[b]{0.5\columnwidth}
        \includegraphics[width=\linewidth]{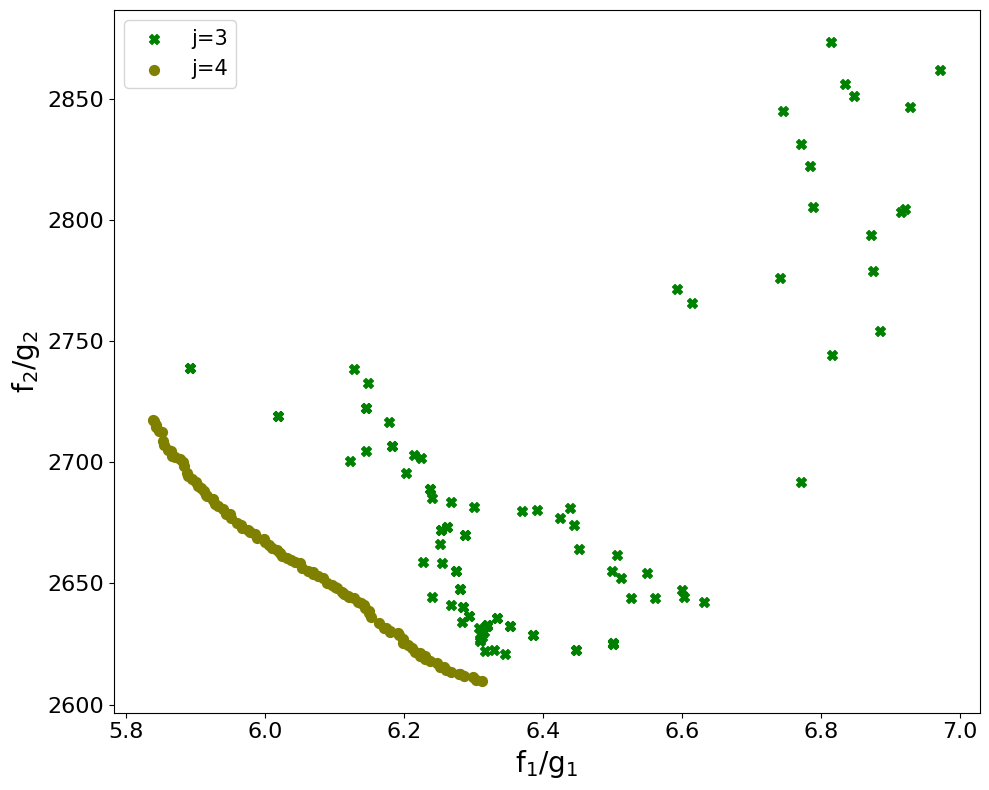}
        \caption{}
        \label{fig:sub7_n20}
    \end{subfigure}%
    \begin{subfigure}[b]{0.5\columnwidth}
        \includegraphics[width=\linewidth]{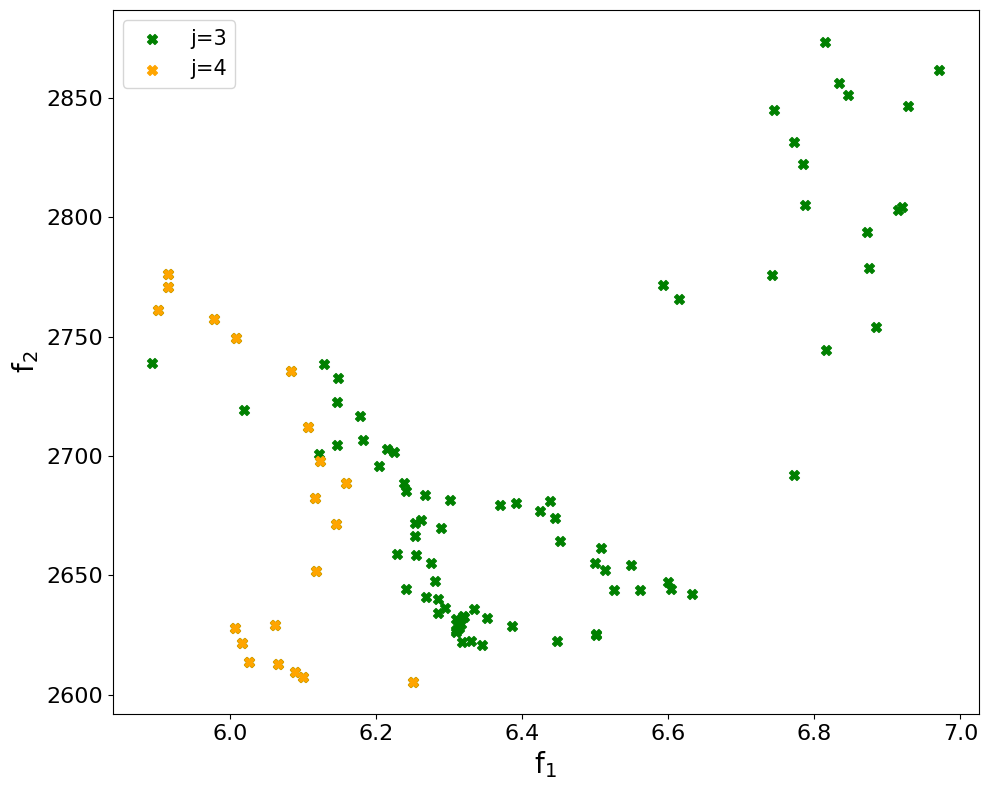}
        \caption{}
        \label{fig:sub8_n20}
    \end{subfigure}
    \begin{subfigure}[b]{0.5\columnwidth}
        \includegraphics[width=\linewidth]{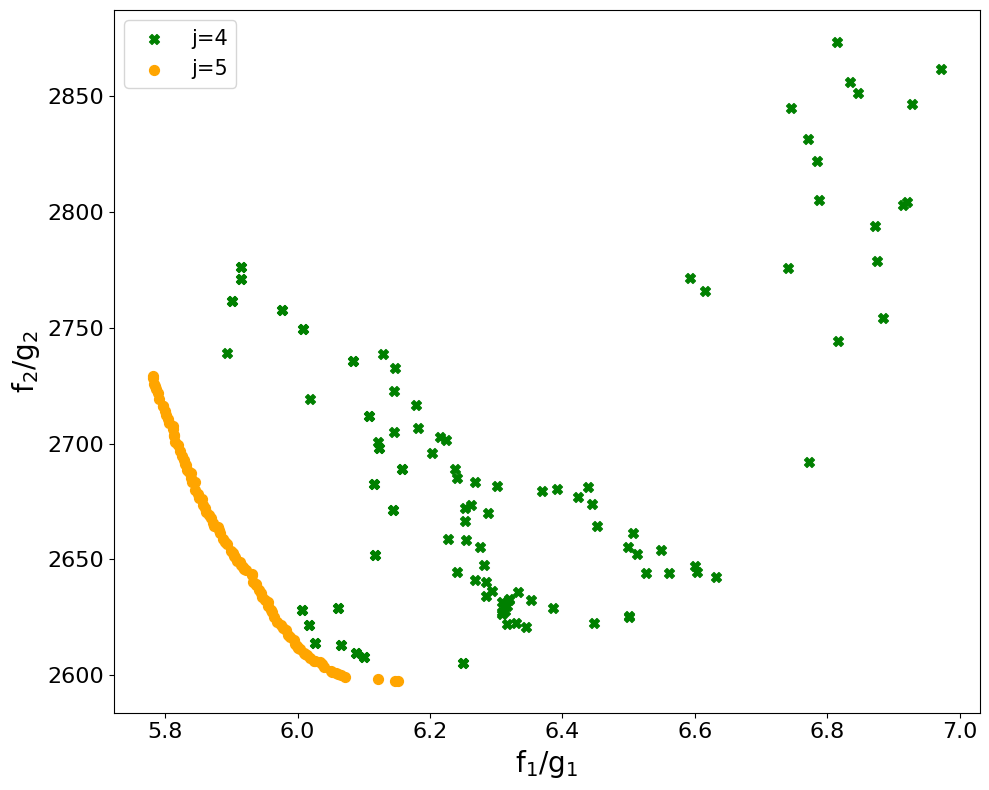}
        \caption{}
        \label{fig:sub9_n20}
    \end{subfigure}%
    \begin{subfigure}[b]{0.5\columnwidth}
        \includegraphics[width=\linewidth]{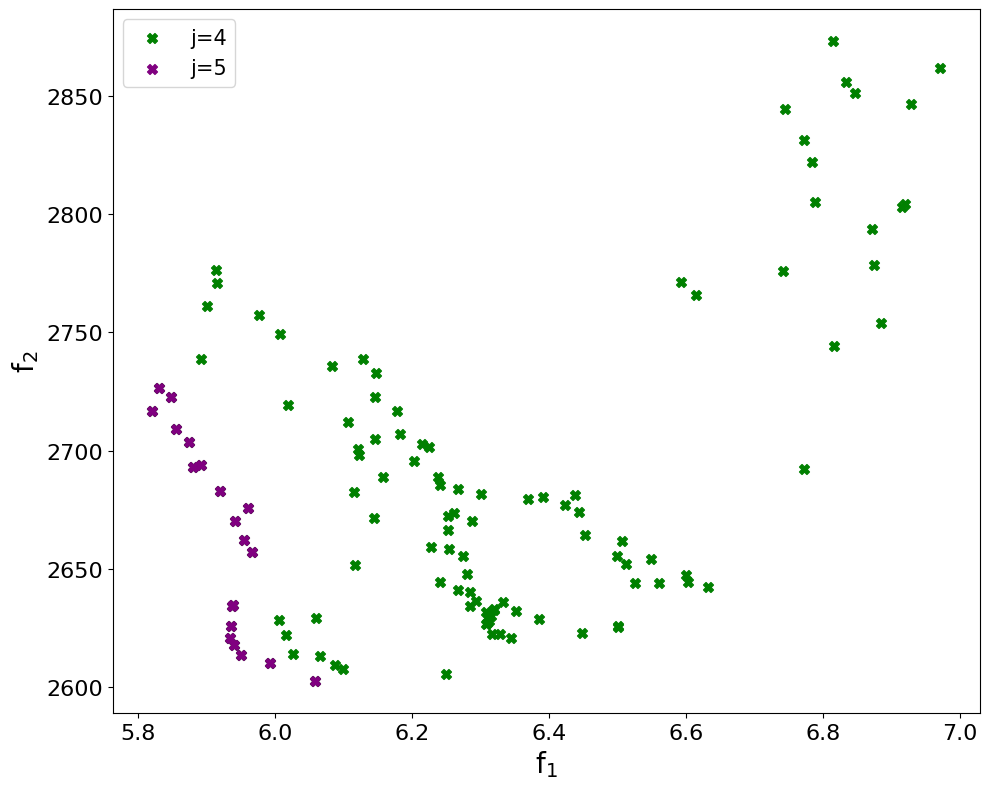}
        \caption{}
        \label{fig:sub10_n20}
    \end{subfigure}%
    \begin{subfigure}[b]{0.5\columnwidth}
        \includegraphics[width=\linewidth]{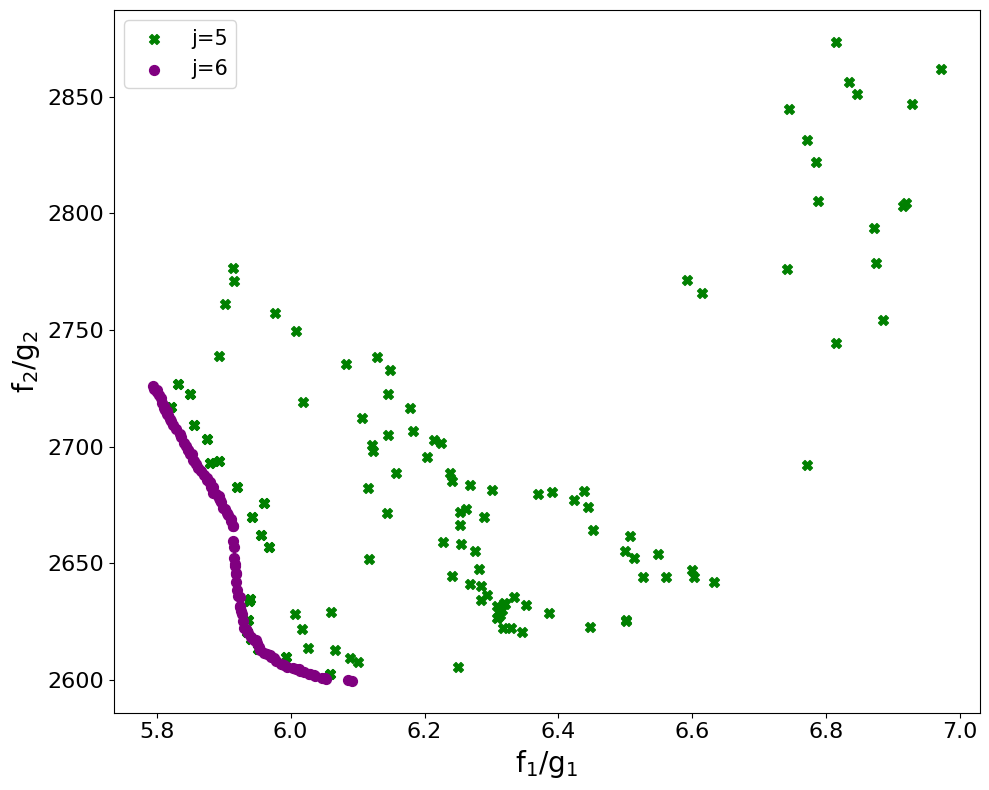}
        \caption{}
        \label{fig:sub11_n20}
    \end{subfigure}%
    \begin{subfigure}[b]{0.5\columnwidth}
        \includegraphics[width=\linewidth]{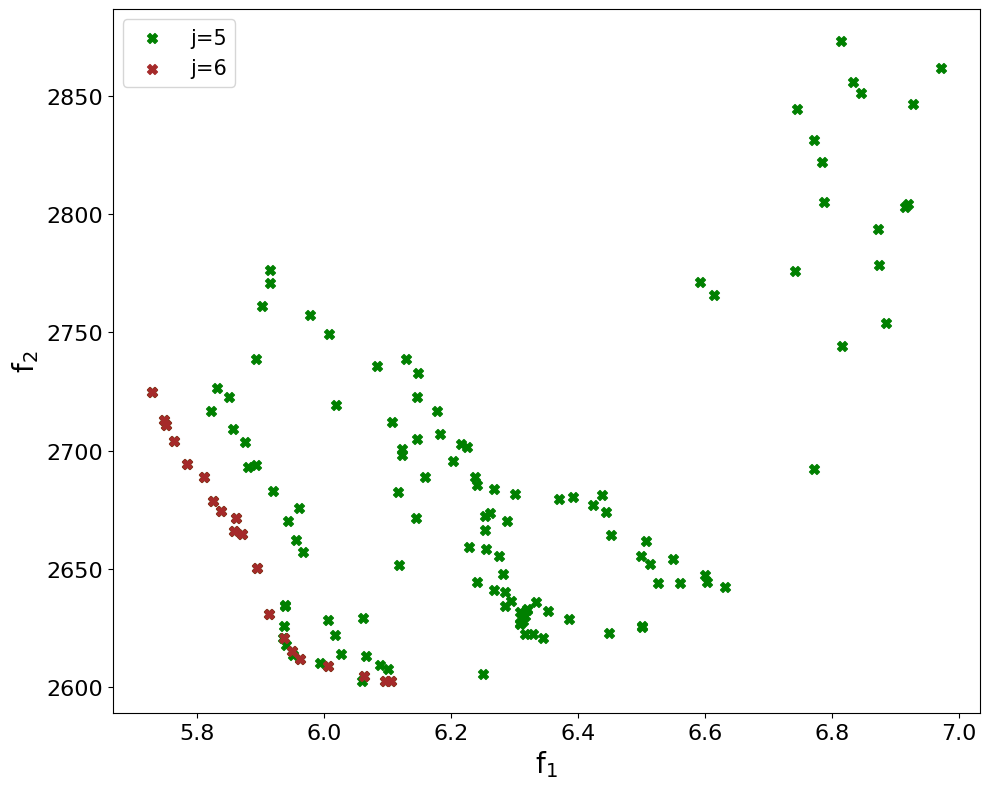}
        \caption{}
        \label{fig:sub12_n20}
    \end{subfigure}
    \begin{subfigure}[b]{0.5\columnwidth}
        \includegraphics[width=\linewidth]{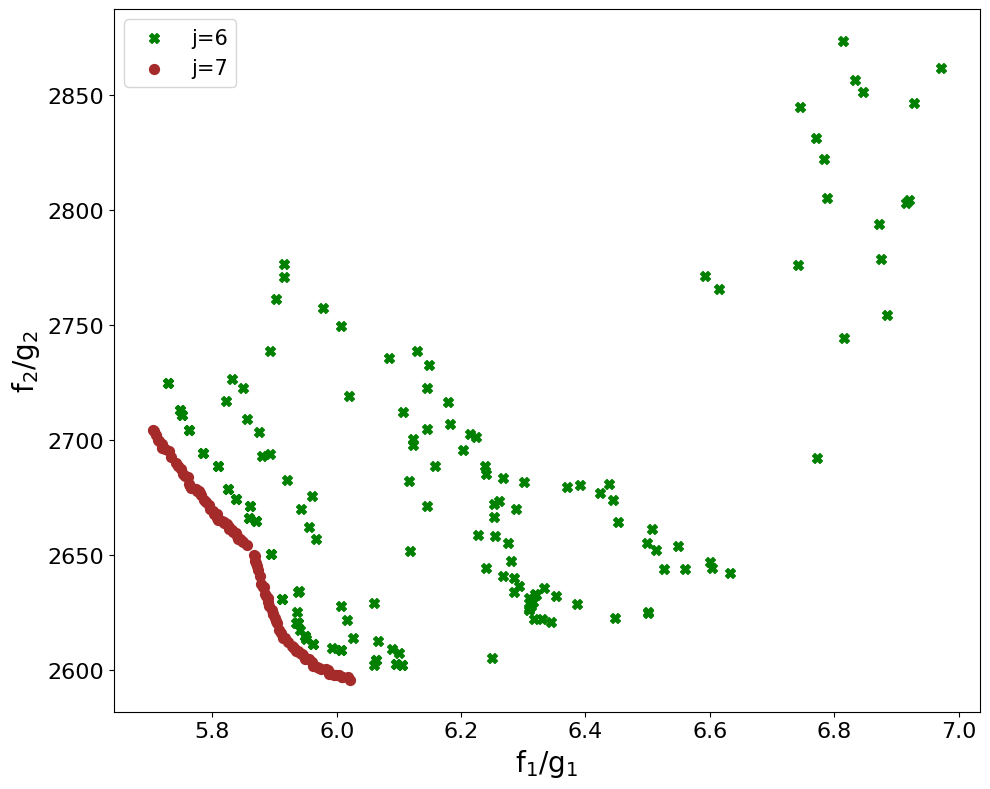}
        \caption{}
        \label{fig:sub13_n20}
    \end{subfigure}%
    \begin{subfigure}[b]{0.5\columnwidth}
        \includegraphics[width=\linewidth]{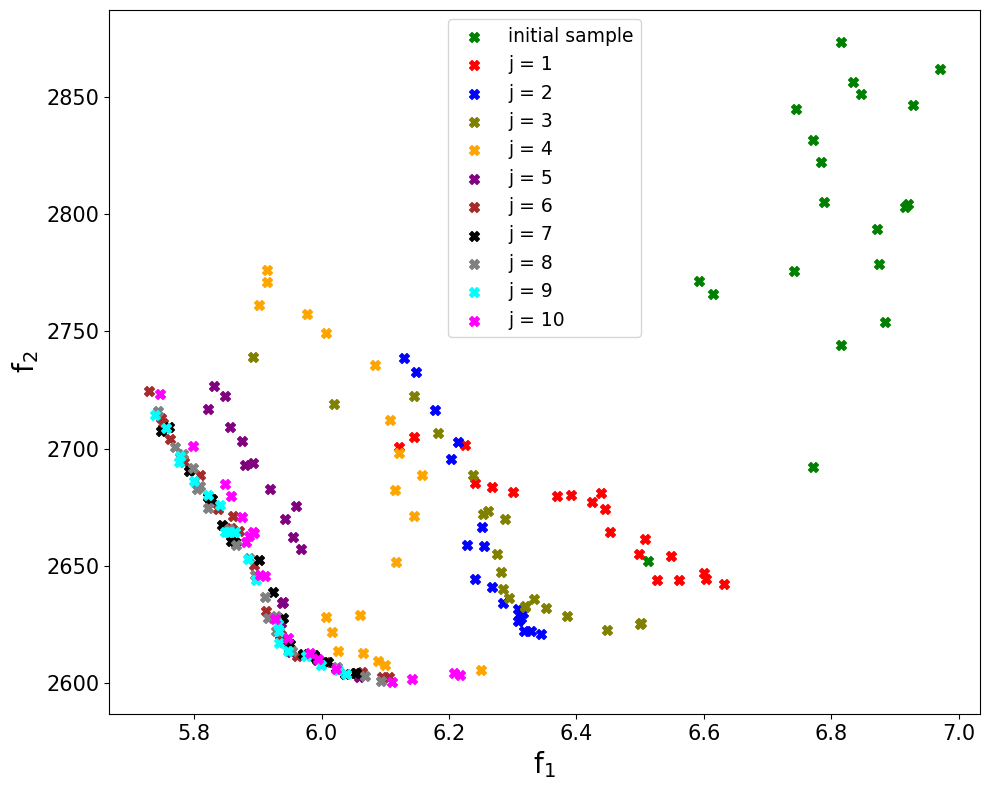}
        \caption{}
        \label{fig:sub14_n20}
    \end{subfigure}
    \begin{subfigure}[b]{0.5\columnwidth}
        \includegraphics[width=\linewidth]{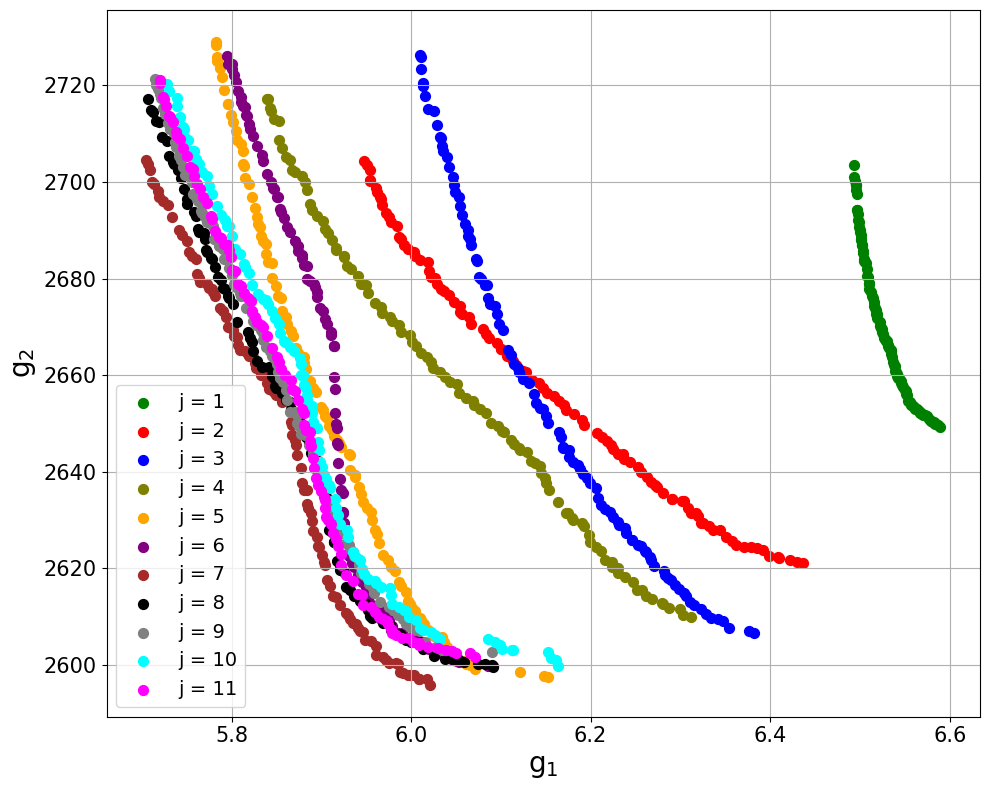}
        \caption{}
        \label{fig:sub15_n20}
    \end{subfigure}%
    \caption{A breakdown of the iterative steps taken to compute $\Pfronthat$ for the MBS. We use $s=20$ as shown in Figure \ref{fig:lhs20}, then always sample $20$ additional points (cf.\ Figure \ref{fig:sampledKM}). We use an ANN surrogate model and NSGA-II. Sub-figures \ref{fig:sub1_n20}--\ref{fig:sub13_n20} illustrate the additional samples at each iteration and the newly computed $\Pfronthat$ until the Hausdorff convergence criterion is met. Sub-figures \ref{fig:sub14_n20} and \ref{fig:sub15_n20} show a summary in different colors of both the sampled points and the Pareto fronts computed in each iteration.}
    \label{fig:NsgaNN20}
\end{figure*}

In our setting with $s=20$ and a convergence threshold of $h_{\mathsf{min}} = 2$ for the Hausdorff distance (a larger threshold yielded a much too early termination of the algorithm), 
$j=11$ iterations were required for convergence. Figure \ref{fig:NsgaNN20} shows a breakdown of the iterative steps. 
Note that expensive model evaluations $f(x)$ are denoted by cross signs, whereas surrogate-based points $g(x)$ are shown as circles . Depending on which quantities are shown, the axes show either $f_i$, $g_i$ or both.
The true objective function values of the $k$-means sampled points are shown in Figure \ref{fig:sub1_n20} in red. In the next iteration, the data points are combined to retrain the surrogate model and compute $\Pfronthat$ with $M=100$ points as illustrated in \ref{fig:sub2_n20}. Since the stopping criterion $h<h_{\mathsf{min}}$ is not satisfied, $20$ new data points are sampled as shown in Figure \ref{fig:sub3_n20} and \ref{fig:sub4_n20}, and combined to compute the $\Pfronthat$ in Figure  \ref{fig:sub5_n20} in blue. The remainder of the sub-figures in \ref{fig:NsgaNN20} shows the iterative process until convergence. The color coding is always of the form 
\begin{itemize}
    \item points $f(x)$ sampled from $\iterate{\Psethat}{j-1}$ (e.g., the blue in \ref{fig:sub4_n20}),
    \item new Pareto front $\iterate{\Pfronthat}{j}$ based on the re-trained surrogate in the same color (Figure \ref{fig:sub4_n20}).
\end{itemize}
Figures \ref{fig:sub14_n20} and \ref{fig:sub15_n20} show summaries of the $S=120$ sampled data points and the corresponding $\iterate{\Pfronthat}{j}$ computed in each iteration. What ultimately matters is that the distance becomes small, both in terms of the collected samples (\ref{fig:sub14_n20}) and the surrogate-based fronts (\ref{fig:sub15_n20}), which nicely demonstrates the convergence criterion $h<h_{\mathsf{min}}$. Note that a smaller value of $h_{\mathsf{min}}$ further imporoves the approximation, but at the cost of way more iterations $j$.

\subsection{NSGA-II with RBF}
\begin{figure}
    \begin{subfigure}{.5\columnwidth}
      \centering
      \includegraphics[width=\columnwidth]{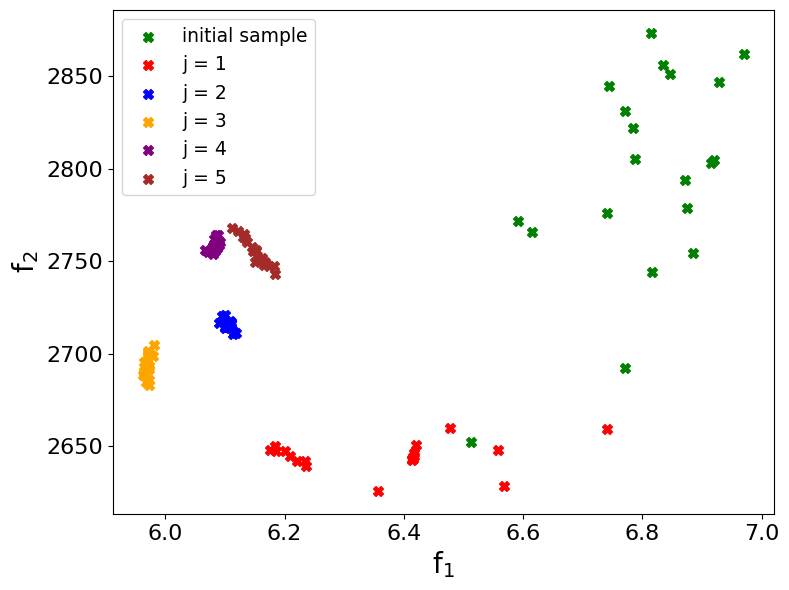}
      \caption{}
      \label{fig:sub1rbf}
    \end{subfigure}%
    \begin{subfigure}{.5\columnwidth}
      \centering
      \includegraphics[width=\columnwidth]{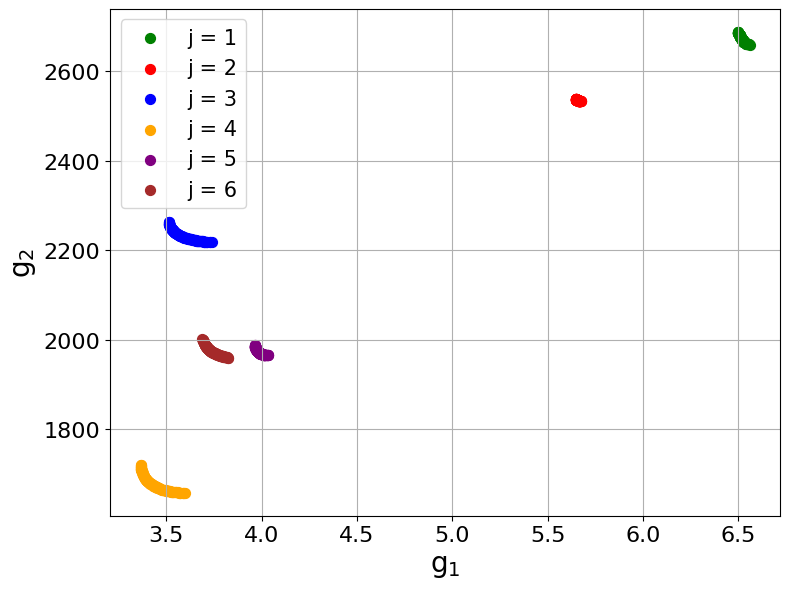}
      \caption{}
      \label{fig:sub2rbf}
    \end{subfigure}
    \caption{(\subref{fig:sub1rbf}) a summary plot of all the different $20$ points computed in each iteration until convergence using the RBF surrogate model and NSGA-II in Algorithm \ref{alg:SAMO}. (\subref{fig:sub2rbf}) illustrates the Pareto fronts  $\Pfronthat$ computed in each iteration.}
    \label{fig:NsgaRBF20}
\end{figure}

For further experimentation, we replace the ANN surrogate model by an RBF model with a Guassian kernel, while the other settings remain unchanged. Different parameter values for the kernel width, i.e., $\sigma\in\{0.1, 0.5, 1.0, 2.0, 5.0\}$, were tested and $\sigma=0.5$ was observed to yield the smallest mean squared error. Using the RBF surrogate model, Algorithm \ref{alg:SAMO} converges after $6$ iterations. 


Figure \ref{fig:NsgaRBF20} shows a summary plots of the sampled data points $f(x)$ and Pareto front $\Pfronthat$.
We observe a worse performance in comparison to the ANN surrogate model, which we believe is due to the fact that RBF models are more sensitive to unevenly distributed samples. For future work, we believe that it will be important to improve sampling strategies to get better results. We believe that this will allow us to leverage the good approximation properties that RBFs have shown in other multi-objective contexts before \cite{Berkemeier2021}.

\subsection{Different sample sizes}
\begin{figure*}
    \centering
    \begin{subfigure}[b]{0.5\columnwidth}
        \includegraphics[width=\linewidth]{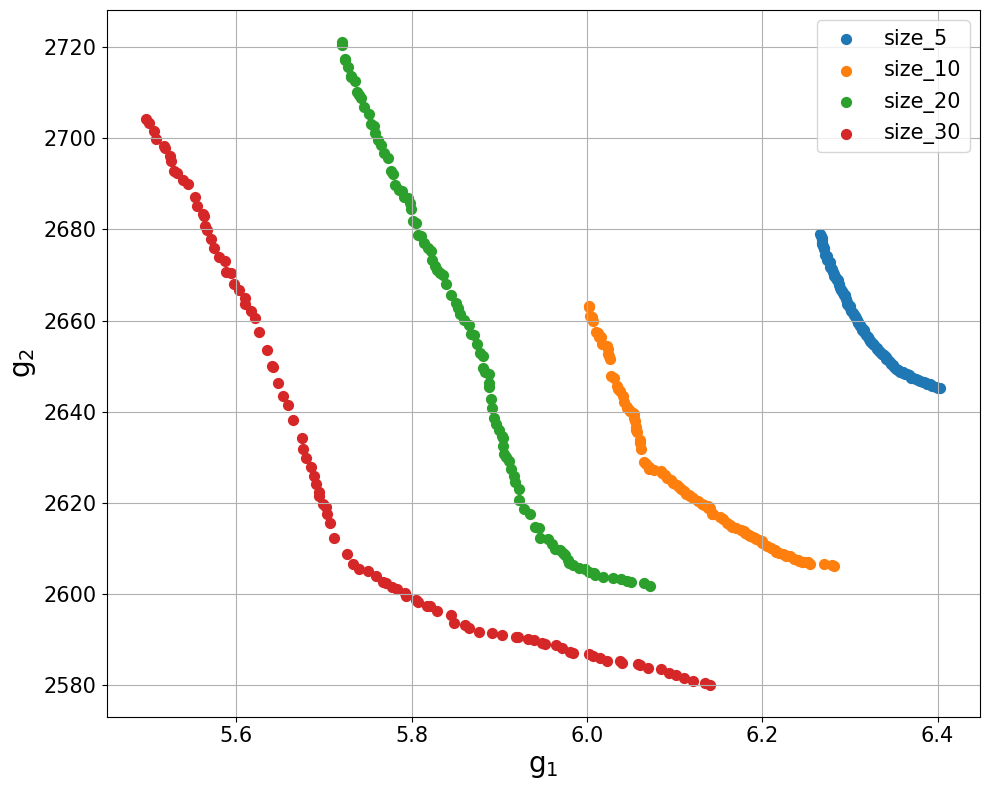}
        \caption{}
        \label{fig:sub1}
    \end{subfigure}%
   \begin{subfigure}[b]{0.5\columnwidth}
        \includegraphics[width=\linewidth]{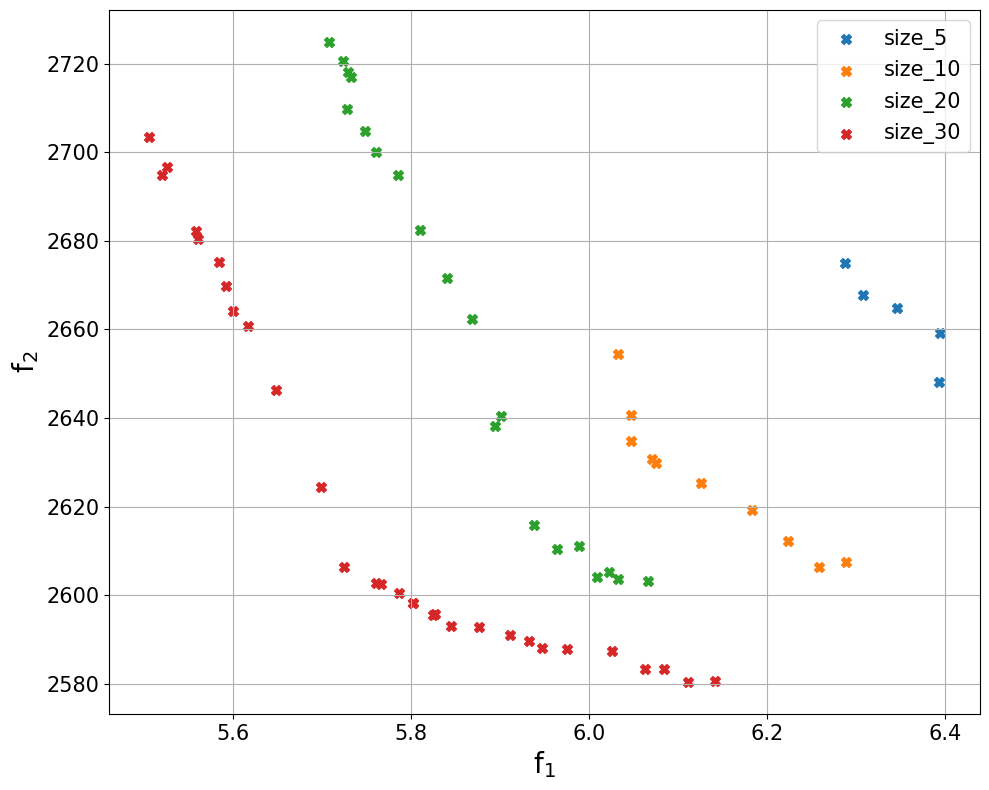}
        \caption{}
        \label{fig:sub2}
    \end{subfigure}%
    \begin{subfigure}[b]{0.5\columnwidth}
        \includegraphics[width=\linewidth]{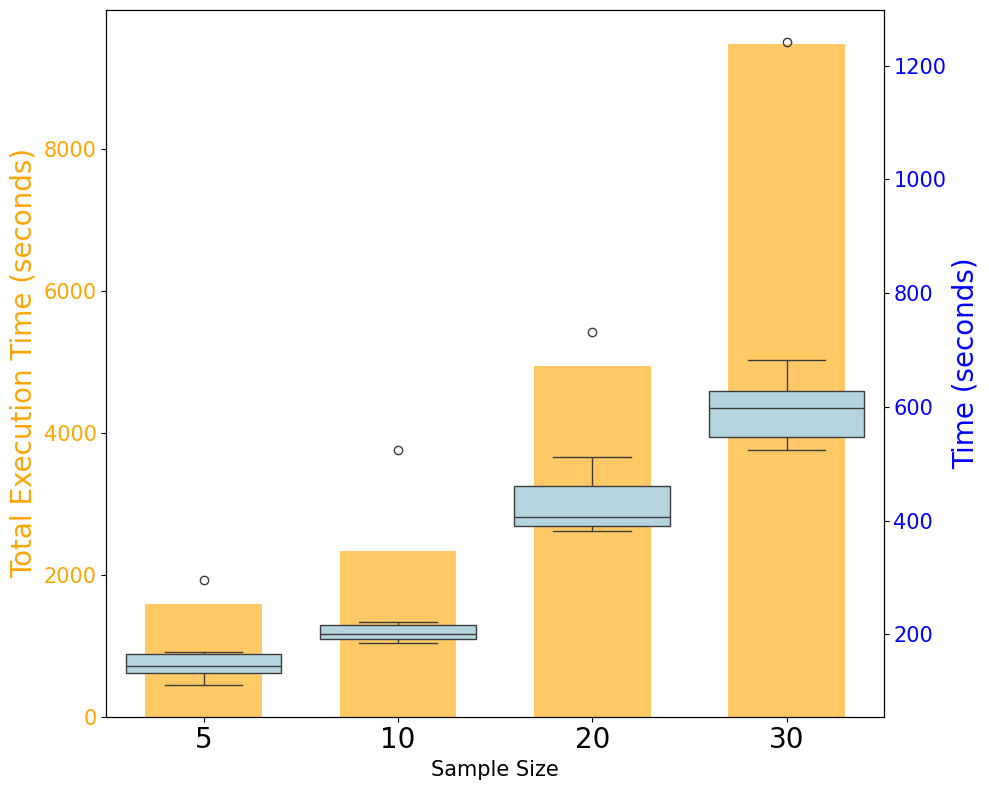}
        \caption{}
        \label{fig:sub3}
    \end{subfigure}%
    \begin{subfigure}[b]{0.5\columnwidth}
        \includegraphics[width=\linewidth]{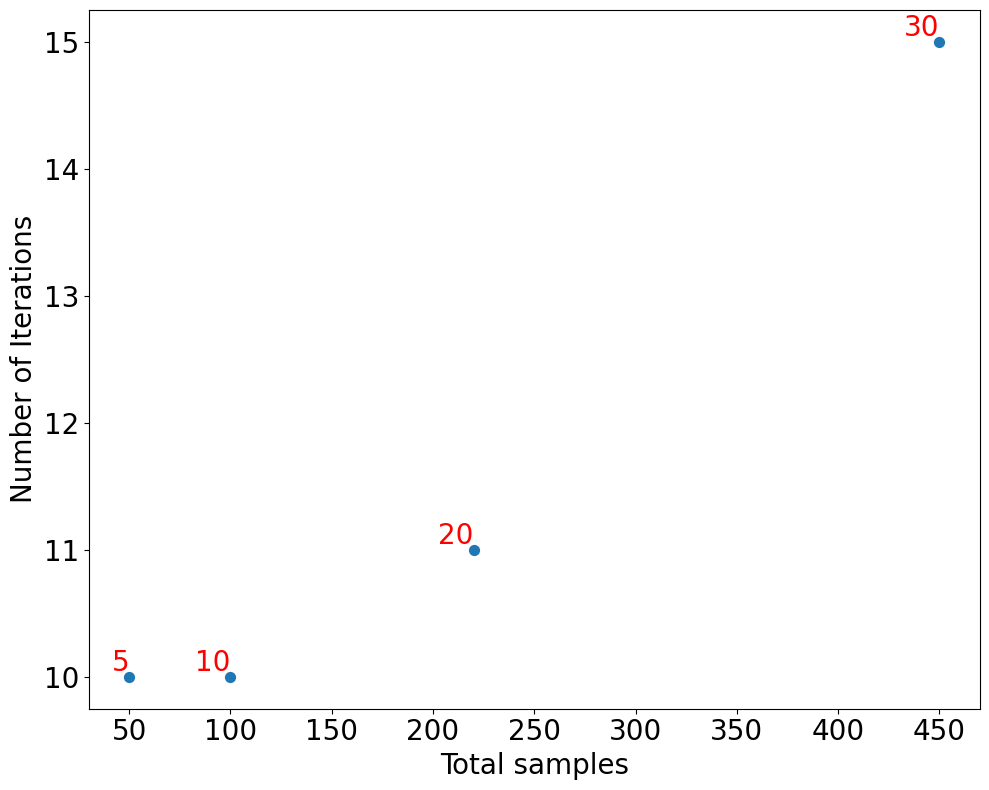}
        \caption{}
        \label{fig:sub4}
    \end{subfigure}   
    \caption{(\subref{fig:sub1}) and (\subref{fig:sub2}) illustrate the final Pareto front $\Pfronthat$ after convergence and the last computed samples of the MBS complex model for the different sample sizes. The computation times (total and per iteration) until convergence is illustrated in (\subref{fig:sub3}) and (\subref{fig:sub4}) respectively. (\subref{fig:sub3}) number of iterations $j$ versus sample size $s$. 
    }
    \label{fig:NsgaNNRBF}
\end{figure*}

Using different sample sizes of $s\in\{5,10,20,30\}$, Algorithm \ref{alg:SAMO} is repeated using both the ANN and RBF surrogate models. The threshold of $h_{\mathsf{min}} = 2$ as a convergence criterion for Algorithm \ref{alg:SAMO} was maintained for all sample sizes. From Figure \ref{fig:NsgaNNRBF} we observe that the larger the sample size, the fewer iteration are needed for convergence, which can be expected. On the other hand, the compute time per iteration, as well as the overall computation time increase. and a better approximated Pareto front. For future research, adaptive sample sizes thus appear to be an interesting approach to allow for quick improvements in the beginning, and for a good final performance in the end.

\section{Conclusion}\label{sec:Conclusion}
We have presented an adaptive surrogate-based framework for solving expensive multi-objective optimization problems. The key ingredient is a Hausdorff measure to study convergence behavior over a back-and-forth procedure between sampling and optimization. We find that---as expected---trade-offs need to be made between approximation accuracy and the number of expensive model evaluations. A significant speedup of several orders was observed while indicating convergence to high-quality solutions.

For future work, we will investigate a larger number of objectives as well as advanced sampling techniques. We here simply used $s$ samples from each iteration, which can result in duplicates or at least clustered sampling over multiple iterations. Moreover, a more efficient surrogate update strategy, where we only use the new samples---in combination with strong regularization or early stopping---to update our model, thus allowing for more iterations and better convergence.

\section*{Acknowledgments}
ACA and SP acknowledge funding by the German Federal Ministry of Education and Research (BMBF) through the AI junior research group ``Multicriteria Machine Learning'' (Grant ID 01$|$S22064).
MBB, MW, WS and SP acknowledge funding by the DFG Priority Programme 2353 ``Daring more intelligence'' (Project ID 501834605).

\bibliographystyle{IEEEtranS}
\bibliography{references}

\end{document}